# Multi-Criteria Scheduling of realistic Flexible Job Shop: A novel approach for integrating Simulation Modeling and Multi-Criteria Decision Making

M.Thenarasu[a*], K. Rameshkumar[a], M. Di Mascolo[b] and S.P. Anbuudayasankar[a]

[a]Department of Mechanical Engineering, Amrita School of Engineering, Coimbatore, Amrita Vishwa Vidyapeetham, India
[b]Univ. Grenoble Alpes, CNRS, Grenoble INP**, G-SCOP, 38000 Grenoble, France
*Corresponding author email id: m_thenarasu@cb.amrita.edu
**Institute of Engineering Univ. Grenoble Alpes

**Abstract**

Increased flexibility in job shops leads to more complexity in decision-making for shop floor engineers. Partial Flexible Job Shop Scheduling (PFJSS) is a subset of Job shop problems and has substantial application in the real world. Priority Dispatching Rules (PDRs) are simple and easy to implement for making quick decisions in real-time. The present study proposes a novel method of integrating Multi-Criteria Decision Making (MCDM) methods and the Discrete Event Simulation (DES) Model to define job priorities in large-scale problems involving multiple criteria. DES approach is employed to model the PFJSS to evaluate Makespan, Flow Time, and Tardiness-based measures considering static and dynamic job arrivals. The proposed approach is implemented in a benchmark problem and large-scale PFJSS. The integration of MCDM methods and simulation models offers the flexibility to choose the parameters that need to govern the ranking of jobs. The solution given by the proposed methods is tested with the best-performing Composite Dispatching Rules (CDR), combining several PDR, which are available in the literature. Proposed MCDM approaches perform well for Makespan, Flow Time, and Tardiness-based measures for large-scale real-world problems. The proposed methodology integrated with the DES model is easy to implement in a real-time shop floor environment.

## 1. Introduction

Production scheduling plays a significant role in the production and service industry as a decision-making operation. It involves a determination of the order of processing a set of tasks on resources over a period of time. Scheduling aims to maximize the efficiency of the operation and reduce costs in the manufacturing and service industries. Owing to the increasing complexity of the production systems, the scheduling techniques have evolved in parallel. Dynamic scheduling is preferred to static scheduling owing to the ability to incorporate uncertainties over time. Dolgui and Proth (2010) presented two algorithms capable of giving dynamic delivery date for assembly and linear production systems. Various types of important scheduling problems include Single Machines, Parallel Machines, Flow Shops, and Job Shops. Gordon et al. (2002) exhaustively reviewed the single and parallel machine models to provide a unified framework for scheduling and due date assignment problem scenarios. A holistic view of the models, algorithms, and properties have been discussed. In Flow Shops, every job follows the same trajectory and gets processed by all the machines in a predefined sequence. In Job Shops, '*n*' jobs are to be processed on '*m*'-machines, each with a unique operating sequence. Each work can only be handled by one machine at a time, and each machine can only do one job at a time. The need for efficient scheduling has recently increased owing to market demand for product quality, flexibility, and timely delivery of the product to the customers.

Owing to the importance of real-time scenarios, the Job Shop Scheduling Problem (JSSP) has been the focus of many researchers. JSSPs are defined such that a set of machines, say $M_j$ ($j = 1, 2, m$) are scheduled to a given set of jobs, say $J_i$ ($i = 1, 2, n$), to optimize the performance measures (Pinedo, 2005). The objective of JSSP is to find an optimal schedule based on criteria such as Makespan, Flow Time, Tardiness, etc. JSSP is one of the hardest Combinatorial Optimization (CO) problems. JSSPs can be categorized as NP-Hard problems wherein an optimal solution cannot be derived even for small-scale instances due to its large solution space. An extension of JSSP is the Flexible Job-Shop Scheduling Program (FJSSP), which is also an NP-Hard problem with more practical applications. Modern production facilities have managed to grow more effectively and efficiently to accommodate the demands of consumers and allow machinery to do multiple operations with flexibility. The FJSSP involves two sub-problems, one that assigns each job to a specific machine (routing) and the other that schedules jobs onto machines over time (scheduling). The FJSSP can be classified into two categories, i.e., Total FJSSP (T-FJSSP) and Partial FJSSP (P-FJSSP). The former allows each operation to be processed by all machines, while the latter has a compatible set of machines for each operation (Kacem et al. 2002).

## 1.1 Flexible Job Shop Scheduling

Solving FJSSP has attracted interest among researchers and shop floor engineers due to their similarities and proximity to real-time problems. The process of manual scheduling in flexible job shops is extremely difficult and time-consuming for industries (Xie et al. 2019) owing to the larger solution space and complexity of using mathematical and other meta-heuristic approaches (Lunardi, and Ronconi 2021). Hillion and Proth (1989) developed an integer linear programming model based on a timed event-graph model to achieve a near-optimal solution to steady state job shops. The study aimed at maximizing the utilization of the bottleneck resources while minimizing their WIP inventories. The JSSP has been solved using exact and heuristic methods with different objectives minimizing Makespan, Flow Time and Tardiness-based measures (Rameshkumar and Rajendran 2018). Chen et al. (1996) employed a branch-and-bound approach to arrive at an optimal solution considering the earliness and tardiness in scheduling. Due to computational complexity, it is difficult to achieve the best solution for large-scale FJSSPs, and it is recommended to obtain near-optimal solutions utilizing heuristic approaches (Zhang et al. 2020). According to Pongchairerks (2021), industries cannot make quick decisions effectively, given the enormous computational time of meta-heuristics. The literature shows that Priority Dispatching Rules (PDRs) and Discrete Event Simulation (DES) approaches are used to solve large-scale FJSSP to obtain approximate solutions with less computational time.

*1.2 Priority Dispatching Rules (PDR)*

Many researchers have proposed PDRs over the years for ranking jobs due to their lower complexity, ease of implementation, lack of tools and technology, and capacity to generate effective results in less time (Zhang and Wang, 2018). PDRs are frequently used to solve real-world scheduling problems due to the large number of jobs and machines, the complexity of the scheduling environment, and the lack of scheduling software (Jayamohan and Rajendran, 2004). There are several scheduling rules in the literature, but no rule, such as Makespan, mean Tardiness, and Flow Time, can perform well across all performance metrics. Three novel dispatching rules such as PT+WINQ/TIS (Process time +Work in next queue/Time in Shop, PT/TIS (Process time + Time in Shop and AT-RPT(Arrival time + Remaining Process Time), were proposed for improving Tardiness and Flow time-based performance measures, with increased efficiency and accuracy in dynamic Job shops and Flow shops (Rajendran and Holthaus 1999). Tay and Ho (2008) classified dispatching rules into Simple Priority Rule (SPR), Weighted Priority Indices (WPI), Composite Dispatching Rule (CDR), and Heuristic Scheduling Rule (HSR). Simple dispatching rules utilize only one criterion: processing time, the number of processes, arrival time, or due date-based parameters. CDR is a combination of two or more SPRs. SPRs use a single rule for a resource/work center to select a job for processing based on specific criteria such as due date, number of operations, arrival time etc. CDRs are a combination of SPRs used on a shop floor. CDRs are implemented by a) different SPRs are utilized at different workstations/resources separately and b) the Composition of several SPRs to evaluate the priorities of jobs waiting in the queue for processing. Composition is achieved by assigning suitable weights for SPRs. They are deploying different SPRs at resources produced better results than applying a single SPR across all the machines for the shop scheduling problems (Barman,1997). The composition of SPRs improved shop performance for large-scale problems, as studied by Holthaus and Rajendran (1997). CDRs can also be used to solve real-time problems with multiple objectives. The challenge lies in getting the right combination of SPRs to design the CDR. It is observed from the literature that a combination of SPRs gave better results than using individual rules. Ozturk et al. (2018) solved multi-objective FJSSP by minimizing Makespan, Mean Flow Time, and Mean Tardiness using CDRs. Results revealed that CDRs outperformed the SPRs for the benchmark problem. The literature review shows that PDRs are used to solve large-scale problems to overcome computational costs without compromising the solution quality. The limitations of individual rules are improved by implementing CDRs.

*1.3 Multi-Criteria Decision Making*

Many priority rules are based on single criteria for job prioritization instead of considering many realistic factors. In any decision-making process, various factors that influence the objectives are identified, and weightage is assigned to each of these factors (Çelen 2014). Multi-Criteria Decision Making (MCDM) is one of the most prominent branches of decision-making. It has grown rapidly and is one of the important tools to arrive at decisions for problems involving multiple criteria (Sun, 2010). An Analytic Hierarchy Process (AHP) based MCDM model was proposed by Kumar et al. (2017) for selecting the foremost PDR in sustainable energy. Seven different PDRs were considered and the authors used AHP to evaluate and choose the best rule. A total of eight performance measures were considered for the evaluation of sequencing rules. The weights for each criterion were assigned using the AHP. The closeness coefficient for each sequencing rule was calculated, and the rules were ranked. Güçdemir and Selim (2018) used MCDM based simulation modelling approach to select the best dispatching rule in a job shop to minimize Makespan, mean Flow-Time, tardy jobs, and variation from customer expectations for different types of problems; many MCDM approaches have been used to make decisions in the field of engineering (Zavadskas et al. 2014). MCDM approaches considered in this study are TOPSIS, EDAS, CP, and WAM. Fuzzy AHP (FAHP) is often used to determine the proper weights for criteria and the MCDM approaches are used to prioritize the jobs. The literature shows that MCDMs should be explored more for shop scheduling applications considering the variety of criteria that affect the system's performance measures.

*1.4 Discrete Event Simulation*

Industries have utilized simulation-based methodologies to model real-time systems and analyze operational performance measurements (Mostafa and Chileshe 2017). Makespan, Flow Time, and Tardiness-related measures are among the well-known performance measures examined by DES. Queue statistics, inventory, resource utilization, and throughput are other performance metrics discussed (Jilcha et al., 2015). To solve a real-time JSSP for decreasing the Makespan, Habib Zahmani and Atmani (2021) suggested a novel method integrating dispatching rules, data mining, and simulation. Kulkarni and Venkateswaran (2017) used DES to perform Simulation-based Optimization (SbO) and compared it to a mathematical technique; they observed that the SbO approach successfully deals with large complex problems with stochastic processing times. From the literature, it is noticed that simulation-based models are useful in solving large-scale real-world problems and have gained popularity as computer capability has advanced. The industrial problem can be represented in a DES model without violating real-world constraints.

*1.5 Research gap*

In this paper, the simulation model has been integrated with the MCDM technique for solving benchmark and large-scale real word problems. Multiple criteria such as due date, process time, setup time, customer priority and number of operations are considered in this paper to prioritize the jobs. Not much research has been done on the composition of the criteria and their influence on the performance measures using MCDMs for a large-scale real-world problem. Discrete event simulation models are developed to capture the dynamic and probabilistic nature of the shop floor. Jobs are prioritized using the well-known MCDM approaches and integrated with a simulation model for evaluating the performance measures of the shop floor by dynamically updating the job priorities. In this paper, we have validated the proposed approach by considering a large-scale (114 jobs & 27 Machines) real-world problem.

Many studies have focused on solving static problems without considering real-world criteria. For arriving at solutions to a real-world problem, criteria such as the number of operations, setup time, predefined sequence of operations, machine flexibility, stochastic activity times, dynamic job arrival patterns and demand, etc., are essential to incorporate in the model. Meta-Heuristics are used to solve FJSSPs to arrive at optimal or near-optimal schedules for a static problem. Most real-world problems involve dynamic and stochastic activities. Simple PDRs and simulation-based approaches are increasingly used to arrive at quality decisions without much computational complexity. SPRs have produced quality solutions for small-size problems involving a single criterion. CDRs combined with multiple SPRs are improving the performance of the shop floor for different measures, namely Makespan, Flow Time, and Tardiness-related measures for problems involving multiple criteria. CDRs have not been well explored to solve scheduling problems involving real-world constraints. MCDM approaches appear to be effective in comparing alternatives and making decisions in various applications involving multiple criteria. However, MCDM approaches were also not been effectively studied for the FJSSPs problems. DES models are helpful in evaluating the performance measures of large-scale problems with stochastic and dynamic activities. Integrating the DES approach with MCDM or CDRs will help the shop floor engineers solve real-world problems involving multiple criteria. The flexibility of the MCDM approach allows industries to assign weights based on the importance of each criterion in their system.

This work is focused on developing PDRs based on MCDM approaches for large-scale PFJS involving static and dynamic job arrivals. The proposed work integrates the MCDM approach with the DES model for PFJSSP containing realistic constraints. Benchmark instances and real-world, large-scale industrial problems are studied using best-known CDRs selected from the literature and proposed MCDM methods. A Fuzzy AHP (FAHP) is applied to assign the weights of realistic criteria considered in this study, and job prioritization is done using MCDM methods. FAHP processes are well-suitable for problems involving bias or preferences whenever human factors assign weights for the criteria. FAHP is an integration of AHP and Fuzzy Logic theory. Unlike in AHP, where the scale of relative importance is represented by the crisp set numerical, FAHP uses fuzzy numbers. DES model is used to compute the Makespan, Flow time, and Tardiness-based performance measures. This study presents four MCDM methods: TOPSIS, EDAS, CP, and PROMETHEE. The best-performing CDRs are employed to prioritize the jobs to minimize Makespan, Mean Flow Time, Mean and Mean and Maximum Tardiness.

*1.6 Objectives*

Based on the detailed literature review, research gaps are identified, and the following objectives are framed:

a. To develop a DES model for PFJSSP to evaluate Makespan, Flow Time and Tardiness-based performance measures
b. To study the performance of CDRs for large-scale PFJSSP
c. To develop hybrid PDRs for job prioritization by integrating the DES model and MCDM methods
d. To implement proposed PDRs based on the MCDM method in a large-scale real-world PFJSSP

The simulation model proposed in this study is modular and flexible, incorporating industry-specific parameters such as the number of resources, number of job variants, static and dynamic job arrival patterns, sequence of the jobs, processing time, and job flexibility. The model proposed in this study is generic and can be used to conduct a 'what if type' analysis considering different scenarios. This analysis will help optimize the available resources and choose the best strategy based on their objective. The DES model developed in this study is integrated with MCDM-based rules and dynamically updates the job ranking as the simulation progresses based on the chosen criteria. The computational complexity of PDRs/CDRs proposed in this study is relatively lesser than heuristics/metaheuristics. The industries can quickly adapt and implement the proposed approach for decision-making. Further, a decision support system for shop floor scheduling can be developed using the proposed methodology based on the industry requirement. The main contribution of this work is summarized as follows:

- Developed Discrete Event Simulation (DES) Model of a large-scale Partial Flexible Shop consisting of Static and Dynamic job arrival scenarios and deterministic and stochastic activity times
- Implemented best-performing Composite Dispatching Rules (CDRs) available in the literature to solve benchmark and large-scale real-world Partial Flexible Shop problem
- Developed Multi-Criteria based Priority Dispatching Rules (PDRs) to solve benchmark and large-scale real-world Partial Flexible Shop problem
- The results obtained in this study could be used as benchmark solutions for the researchers working on developing new PDRs for PFJS.

The remainder of the paper is structured as follows. The mathematical formulation of the FJSSP is given in section 2. The methodology for developing the simulation model is presented in section 3. The proposed MCDM approaches and their implementation methodology for ranking jobs based on real-world criteria are discussed in section 4. Performance evaluation of PFJS using the best-known CDRs available in the literature and proposed MCDMs are presented in section 5. The proposed MCDMs are implemented in a real-world PFJS. Their performances are compared by evaluating well-known performance metrics, namely Makespan, Flow Time and Tardiness-based measures presented in section 6. Conclusions and future directions of the research are presented in section 7.

## 2. Mathematical Formulation of Partial Flexible Job Shop Scheduling Problem

A PFJSSP is a scenario that handles *n* jobs using *m* machines, with each operation eligible to be performed on a set of alternate machines. The jobs are processed based on a predefined sequence using any of the feasible machines for each operation. FJSSP involves assigning operations to a machine and finding the processing order of jobs on each machine. The key objective is to improve the performance indicators Makespan, Maximum Tardiness, and total Flow Time. The mathematical model of Partial Flexible Job Shop Scheduling (PFJSS) was derived by De Araujo and Previero, 2019.

Notations used in the mathematical model:

| | |
|---|---|
| $M_{max}$ | : Makespan |
| $F_{mean}$ | : Mean Flow time |
| $T_{mean}$ | : Mean Tardiness |
| $T_{max}$ | : Maximum Tardiness |
| $n$ | : Number of jobs; set of jobs: $J = \{J_1, J_2, \ldots, J_n\}$ |
| $m$ | : Number of machines; set of machines: $M = \{M_1, M_2, \ldots, M_m\}$ |
| $p$ | : Number of operations; set of operations: $O = \{O_1, O_2, \ldots, O_p\}$ |
| $i$ | : Job index: $i=1,\ldots n$ |
| $i'$ | : Preceding job index: $i'=1,\ldots,n$ |
| $j$ | : Machine index: $j=1,\ldots m$ |
| $k$ | : Operation index: $k=1,\ldots p$ |
| $f_j$ | : Flexibility level of machine '$j$' |
| $L$ | : Upper bound on the Makespan of the best solution |
| $P$ | : Set of operation $k$ and its preceding operation $k'$ |

| $M$ | : Set of all machines for operation $k$ |
|---|---|
| $E$ | : Set of all compatible machines |
| $C_i$ | : Completion time for job $i$ |
| $T_i$ | : Tardiness of job $i$ |
| $F_i$ | : Flow time of job $i$ |
| $D_i$ | : Due date of job $i$ |
| $Pt_i$ | : Total Processing time of job $i$ |
| $Pt_{i'}$ | : Total Processing time of preceding job $i'$ |
| $S_i$ | : Start time of an operation of job $i$ |
| $S_{i'}$ | : Start time of a preceding operation of job $i'$ |
| $E_{ik}$ | : Number of compatible machines for operation $k$ of job $i$ |
| $Pt_{ijk}$ | : Processing time of operation $k$ of job $i$ on machine $j$ |
| $x_{ijk}$ | : Binary variable with value 1 if operation $k$ of job $i$ is assigned on machine $j$; 0 otherwise |
| $y_{kk'}$ | : Binary variable with value 1 if operation $k$ precedes the operation $k'$ on a given machine. Else, holds value 0 or 1 |

Minimization of

$$M_{max} = \text{Min}(\text{Max}(C_i)) \; ; \; F_{mean} = (\sum F_i)/n; \; T_{mean} = \sum_{i=1}^{n}\frac{T_i}{n} \; ; T_{max} = \max_i T_i \quad \text{---------------(1)}$$

The objective function is to minimize Makespan, Mean Flow time and Tardiness-related measures described in equation (1)

Subject to:

$$\sum_{j \in m} x_{ijk} = 1, \forall i \in J \; \forall k \in O \quad \text{-------------------------------------------------------------(2)}$$

Where $x_{i,j,k}$ represents operation $k$ of job $i$, on machine $j$. Sum of probability of an operation of a job across machines is always 1. This is because an operation of a job can be processed by one and only one of the set of all feasible machines.

$$Pt_i = \sum_{k \in p, j \in m} Pt_{ijk} \cdot x_{ijk}, \forall i \in J \; \forall k \in O \quad \text{------------------------------------------(3)}$$

Where $Pt_{ijk}$ represents processing time of operation $k$ of job $i$, on machine $j$. Given the constraint one, and that the term $x_{ijk}$ is binary and takes value 1 only in one instance for any job-operation-machine combination, it is clear that that process time can only be a single term and is equal to $Pt_i$ or $Pt_{ijk}$ for operation $k$ in the feasible machine $j$ that performs the particular job $i$

$$S_i + Pt_i \leq M_{Max}, \forall i \in J \quad \text{(4)}$$

Where $S_i$ and $Pt_i$ are the start time and process time of job $i$

The maximum possible value that a process can take after the start of an operation for a job in the system is always such that the time elapsed is less than or equal to makespan.

$$y_{kk'} + y_{k'k} \geq x_{ijk} + x_{i'jk} - 1,$$
$$\forall k \in M, (k, k') \in E_{ik} \quad \text{(5)}$$

Precedence checked for all machines. There can be only one possible assignment to any process of job in the system. Also, there can be maximum precedence of only one operation when comparing between two operations of a job.

$$S_{i'} + Pt_{i'} - L(1 - y_{kk'}) \leq S_i, \forall (k, k') \in E_{ik} \quad \text{(6)}$$

The start time of all preceding operations ($S_{i'}$) and total Processing time of preceding job $i'$ is always less than or equal to the total time elapsed before next operation. Where $S_i$ represent the start time of job $i$ and constant $L$ is an upper bound on the best solution's Makespan.

$$S_{i'} + Pt_{i'} \leq S_i, \forall (k, k') \in P \quad \text{(7)}$$

The start time of a preceding operation ($S_{i'}$) and the processing time is always less than or equal to the start time of the next operation.

$$\sum_{i \in O} x_{ijk} \leq f_j \cdot r, \forall k \in M \quad \text{(8)}$$

Includes the flexibility factor of the job shop. $r=1$ when an operation can be performed by all the machines; $0<r<1$ as some machines cannot perform some operations.

$$S_i \text{ and } Pt_i \geq 0, \forall i \in O \quad \text{(9)}$$

The start time and processing time of all the operations of all jobs are non-negative

$$x_{ijk} \in \{0,1\}, \forall k \in M, \forall i \in O_k \quad \text{(10)}$$

Binary variable takes the value 0 when an operation is not performed on a machine and takes value 1 when an operation is performed on a machine.

$$y_{kk'} \in \{0,1\}, \forall (k, k') \in E \quad \text{(11)}$$

Binary variable that indicates the precedence of operation $k$ over operation $k'$. Takes the value 1 if $k$ precedes $k'$; 0 otherwise

## 3. Simulation Model Development of Partial Flexible Job Shop

Arena, a DES software package, has been used to develop the simulation model of PFJSSP. The simulation model developed in this study has been used to a) Simulate the job shop production scenario, b) Implement the decision given by MCDMs to evaluate the performance measures, c) Scheduling of jobs and d) Simulate dynamic and probabilistic activities of the shop floor. A discrete event simulation software 'Arena' has been used to model the benchmark and real-world problem. The model is flexible enough to accommodate variations in the Number of jobs, Number of machines, activity times and job priorities. The simulation model proposed in this study can be easily integrated with the MCDM approach to prioritize the jobs. Based on the priority, a simulation will be carried out, and the performance of the shop floor will be measured. The simulation methodology of PFJS is depicted in Figure 1.

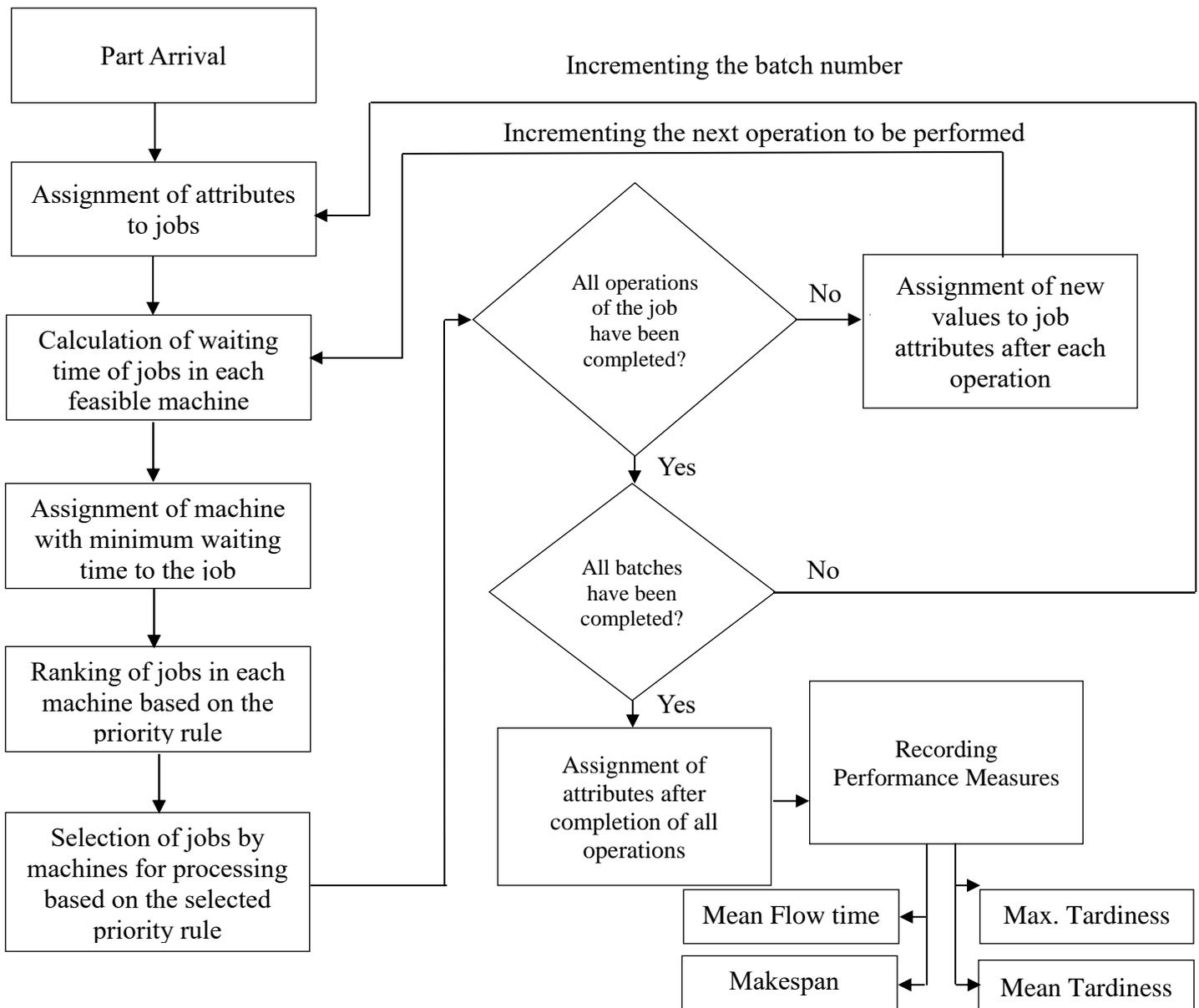

Figure 1 Caption: The trajectory of jobs during the simulation

Figure 1 Alt Text: Functional Flow chart indicating the various steps undergone by incoming jobs, decisions, and collection of performance measures.

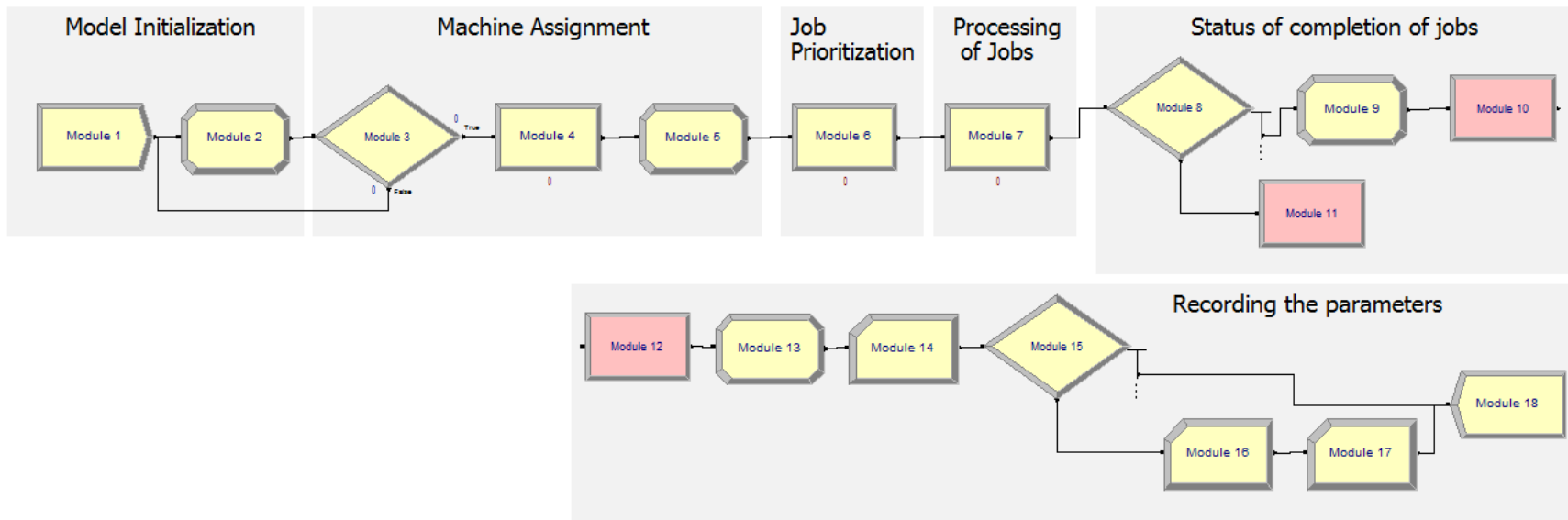

Key for Modules
1. Define the arrival pattern of incoming jobs
2. Assignment of attributes, start time, and index for iteration
3. List and compare the queue lengths for each eligible machine. Is the queue time in a current machine the least?
    a. If Yes: Route to module 4
    b. If No: Determination of queue time on each machine and finding the shortest queue time
4. Machine assignment based on queue
5. Assignment of process time
6. Job selection based on priority rule chosen
7. Job processing
8. Check if all the operations are completed?
9. Increment the attributes to next operation
10. Route to exit station
11. Route to Machines (Job Prioritization)
12. Exit Station
13. Computing completing time of Jobs
14. Recording Makespan
15. Checking the completion time of jobs
16. Counting the late jobs
17. Computing the Tardiness
18. Dispose (Leave the system)

Figure 2 Caption: Simulation model of Partial Flexible

Figure 2 Alt Text: The Snapshot of the DES model drawn using Arena, with descriptions to each module in the Simulation model

A machine loading rule is triggered for choosing a machine for operating when a new job has arrived, when an existing operation is finished, or when its successive operation becomes "operation ready". During this process, a sequencing rule is triggered for choosing the job to be processed subsequently. In this way, each job is mapped to one of the feasible machines based on the least waiting time machine loading rule. Once the sequencing criterion is met, routing action is executed. This signifies that there's an active interplay between both the routing and sequencing processes. For every job operation, the queue time in a set of all feasible machines is iterated, and the assignment is done to the machine with the least queue time. Similarly, iteration happens for all the jobs waiting in each machine to be sent for processing.

An overview of the simulation model, including all the stages, is shown in Figure 2. The model allows flexibility to incorporate changes in job arrival patterns, cycle time, setup time, number of jobs and operations, and priority rules as required. The scalability allows the users to use the model for different problem instances. A statistical tool in Arena DES software generates the distributions for setup time. The model variables, job attributes, and ranks are updated after every operation. This allows the ranking of jobs dynamically as the simulation progresses.

Table 1: Benchmark data and real-world problem instances used for simulation experiments

| FJSSP Instances | | | | | | |
|---|---|---|---|---|---|---|
| Instances | $(n \times m \times k)$** | NOP | Flexibility Factor | Processing time (T.U) | Extension (T.U) | |
| | | | | | Due Date | Setup time |
| $MK$1 | 10×6×55 | 5 to 7 | 3 | 1 to 7 | 16-42 | 0.35-1.14 |
| $MK$2 | 10×6×58 | 5 to 7 | 6 | 1 to 7 | 20-48 | 0.54-1.68 |
| $MK$3 | 15×8×150 | 10 to 10 | 5 | 1 to 20 | 42-84 | 6.24-10.23 |
| $MK$4 | 15×8×90 | 3 to 10 | 3 | 1 to 10 | 21-93 | 3.42-10.55 |
| $MK$5 | 15×4×106 | 5 to 10 | 2 | 5 to 10 | 53-107 | 8.71-17.74 |
| $MK$6 | 10×15×150 | 15 to 15 | 5 | 1 to 10 | 64-116 | 5.92-12.48 |
| $MK$7 | 20×5×100 | 5 to 5 | 5 | 1 to 20 | 80-130 | 3.62-13.27 |
| $MK$8 | 20×10×225 | 5 to 10 | 2 | 5 to 20 | 76-148 | 6.45-16.34 |
| $MK$9 | 20×10×240 | 10 to 15 | 5 | 5 to 20 | 80-154 | 5.78-18.62 |
| $MK$10 | 20×15×240 | 10 to 15 | 5 | 5 to 20 | 84-162 | 8.74-20.45 |
| Real-World FJS problem* | 114×28×245 | 1 to 14 | 4 | 2 to 32 | 120-360 | 0.12-0.48 |

*Note:* ** $n$ – *Number of Jobs,* $m$- *No of Machines and* $k$- *Number of Operations,* **NOP**- *Number. of operations.* **T.U**- *Time Units.* *Time units for real-world FJS problems are in hours*

The simulation is run for ten benchmark problems (*MK*1 to *MK*10), selected from the literature and illustrating FJSSP with various sizes (Brandimarte,1993). The number of jobs, machines, and operations vary from 10 to 20, 6 to 15, and 55 to 240, respectively, and real-time large-scale PFJSSP involving 114 jobs and 28 machines details are given in Table 1. The simulation model developed in this study can be extended as a decision support system for shop floor decision-makers. One static and four dynamic job arrival patterns are simulated to study the effectiveness of the proposed rules. In the static case, it is assumed that all the jobs are available at the beginning of the process.

Table 2: Dynamic Arrival Pattern Scenario

| Dynamic Arrival Pattern | Description | Equation for arrival of jobs |
|---|---|---|
| Equal Interval arrival | Jobs arrive at equal time intervals. | $\sum_{j=1}^{n} j \times Lp$ |
| An increasing rate of arrival | Jobs arrive with an increasing rate of arrival within a defined time frame. | $\sum_{j=1}^{n} t \times \sqrt{\frac{j}{n}}$ |
| Decreasing rate of arrival | Jobs arrive with a decreasing rate of arrival within a defined time frame. | $\sum_{j=1}^{n} t \times \left(\frac{j}{n}\right)^3$ |
| Random arrival | Job arrivals follow Uniform Distribution. | Unif [0, *t*] |

*Note:'**t**' – Total time interval for jobs to arrive, '**n**' – Total number of jobs.*

The dynamic arrival patterns have been used as defined by Sels et al. (2012) and are shown in Table 2. MCDM approaches such as PROMETHEE, CP, EDAS, and TOPSIS are employed to rank and prioritize the jobs considering the benchmark problems and the case instance. The results were compared with the best-known CDRs available in the literature. Table 3 shows the details of CDRs and proposed heuristics to solve the FJSSPs. The composite dispatching rule for multiple objective dynamic scheduling combines different strengths of simple rules for addressing respective objectives. Development of such CDR involves three steps firstly, identifying the elementary practices, then combining them into a single rule, and then estimating the proper values of weights to ensure maximal production performance. Rules C1- C4 are based on weighted priority index rules. These rules are a linear combination of SPRs described with computed weights (Sels et al. 2012). The weight is determined based on the specific business domain and the importance of the job. Similarly, rules C5-C7 are retrieved based on Automatically Defined Functions (ADF). This approach is made by using gene expression programming (GEP) and simulation established by Ozturk et al. (2018). Finally, proposed rules (C9-C12) are developed using a combination of various MCDM methods and the DES model.

Table 3: Composite Dispatching Rules (CDRs) and MCDMs-based priority rules

| Rules | CDRs | Reference |
|---|---|---|
| C1 | 2PT+LWR+FDD | (Sels et al. 2012) |
| C2 | 2PT+LWR+Slack | |
| C3 | SPT+LWR+Slack | |
| C4 | 2PT+LWR+EDD | |
| C5 | (7*LTWC)+(11*SPT)+12*(LNOP+AT) | (Ozturk et al. 2018) |
| C6 | LTWC/(3+LNOP-LRNOP) | |
| C7 | ODD+RT | |
| C8 | [EDD+[(LRNOP+LTWC)/(LRWC-LTWC)]*LNOP]*LRNOP | |
| C9 | FAHP + TOPSIS | Proposed MCDM approach |
| C10 | FAHP + CP | |
| C11 | FAHP + EDAS | |
| C12 | FAHP + PROMETHEE | |

*Note: **PT**- Processing Time, **LWR**- Least Work Remaining, **FDD** – Flow Due Date, **EDD**- Earliest Due Date, **LTWC**- Least Total work content, **SPT**- Shortest Processing Time, **LNOP**- Least No of Operations, **LRWC**- Least Remaining work content, **EDD**- Earliest Due Date, **LRNOP**- Least Remaining number of operations, **PT** – Process Time, **RT**- Remaining Time, **ODD** – Operational Due Date*

## 4 Development of MCDM-based PDR

In this work, MCDM-based PDRs are developed and integrated with the simulation model for arriving job priorities in a PFJS. The PFJSSP identifies factors that influence the system's performance and assigns the corresponding criterion weights. Figure 3 shows the methodology adopted to create a ranking of jobs using the MCDM approach. Delphi method (Dalkey and Helmer, 1963), a powerful technique to assess multiple criteria/factors that influence a decision-making process, is utilized in this study for ranking the criteria. In the Delphi approach, a questionnaire is prepared and circulated to the shop floor engineers and managers to get their feedback. The feedback from shop floor Engineers and Managers is consolidated by FAHP to arrive at the weightage for the criterion. The criteria for ranking jobs are due date, setup time, number of operations, processing time, and slack time per remaining operation. These parameters are selected based on literature review, and data collection carried out in the industry, and discussion with the subject experts. Assigning weights to the criteria is probably the most important part of MCDM, as it allows different views and their impact on ranking alternatives to be explicitly expressed.

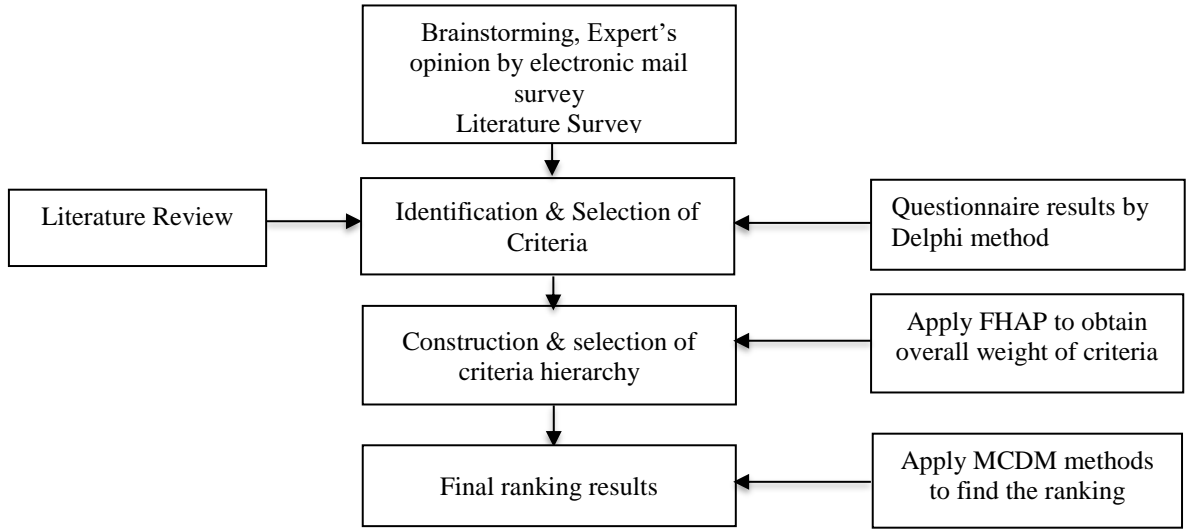

Figure 3 Caption: Methodology of MCDM approach

Figure 3 Alt Text: The overall flow involved in ranking the jobs using MCDM approaches

FAHP processes are suitable for problems involving bias or preferences whenever human factors are involved in assigning weights for the criteria. Details of the computational procedure to arrive at criterion weights are available in Thenarasu et al. (2022). Fuzzy weight details and Crisp weights are shown in Table 4.

Table 4: FAHP Criterion's Weight

| Criteria | Fuzzy Weights | | | Crisp Weights | Crisp Normal Weights |
|---|---|---|---|---|---|
| | $l$ | $m$ | $u$ | | |
| $C_r1$ | 0.07 | 0.10 | 0.16 | 0.34 | 0.10 |
| $C_r2$ | 0.17 | 0.26 | 0.43 | 0.85 | 0.26 |
| $C_r3$ | 0.09 | 0.15 | 0.23 | 0.47 | 0.15 |
| $C_r4$ | 0.03 | 0.04 | 0.06 | 0.12 | 0.04 |
| $C_r5$ | 0.25 | 0.46 | 0.77 | 1.47 | 0.45 |
| Sum of Crisp Weight and Crisp Normal Weight | | | | 3.26 | 1.00 |

*Note: $C_r1$- Process time, $C_r2$- Due Date, $C_r3$- Number of operations, $C_r4$-Setup-time, $C_r5$- Slack time Remaining per Operation, (l m u) – Triangular Fuzzy number, l- Smallest possible value, m-Midpoint value u- Largest value*

*4.1 MCDM Methods for Prioritization of Jobs*

This section describes the job ranking methods based on MCDM methods such as TOPSIS, EDAS, CP and PROMETHEE. The MCDM methods have been adopted in this study to prioritize the jobs in an FJSSP to minimize Makespan, Flow Time, and Tardiness performance measures. Considering the robustness of TOPSIS, it is employed to find decision alternatives for the given criteria based on the distance from an ideal solution (Hwang and Yoon 1981). The input for the TOPSIS approach includes criteria selected for ranking using Delphi and their weights computed using FAHP. Distances from the Positive Ideal Solution (PIS) and the Negative Ideal Solution (NIS) are used to rank the jobs. Upon completion of an operation in a machine, the ranks are updated dynamically. An external spreadsheet is used to compute the job priorities and integrated with the DES model to evaluate the performance measures. EDAS is a relatively new MCDM technique used to decide in situations where there is a conflict among the criteria (Ghorabaee et al. 2015). The ideal alternate is determined based on the distance from the mean solution. Positive Distance from Average (PDA) and Negative Distance from Average (NDA) are used to represent the distance between the mean and the alternate solution. Better-performing alternates are those with higher PDAs and lower NDA. The ranking is done in decreasing order of the appraisal score. The CP is used to attain a solution closest to an ideal solution for a particular distance measure. The distance between these points is the basis for decision-making and transformed evaluation matrix, distance $L_p$– metric represents the Euclidean distance measure. It is computed for all jobs based on the ideal, minimum, and maximum values of each criterion alongside the actual and normalized weights of the selected criteria. The ranking of the alternatives is determined by considering the lowest '$Lp$' metric value from the ideal solution (Diaz-Balteiro et al. 2011). PROMETHEE can be easily implemented due to its computational efficiency in prioritizing jobs in a large-scale problem (Munda 2005). The evaluation matrix lists all the production parameters for developing the PDR. In almost all the practical scenarios, not all parameters would be of the same trend, i.e. some parameters need to be minimized, whereas the others need to be maximized. A good solution may be a balance of all the parameters under consideration based on the business requirements. Since all data are in various trends, it is always essential to normalize the available data to obtain useful and meaningful interpretations. For every job, the preference index and net Flow value are computed for ranking the jobs. Note that the pseudocodes of TOPSIS, EDAS, CP and PROMETHEE method are given in Appendix 1. The use of pseudocode is to rank the jobs based on the MCDM methods chosen. These ranks serve as the input for simulation. This code is run dynamically during the simulation after the completion of each operation. The code is applicable whenever an MCDM is chosen as the preferred PDR. The program runs every time an operation is completed and ranks the jobs based on updated attributes. A similar manner for all job instances ranking is computed using the MCDM approach. The ranking of jobs is calculated dynamically in the external spreadsheet and updated to the DES model. The computation steps of the PROMETHEE method are given below, which is used for prioritizing.

### *4.1.1 Computation of job ranking using the PROMETHEE method*

The step-by-step approach to implementing the PROMETHEE is discussed subsequently. Firstly, the criteria are compared pairwise, followed by the determination of preference function, computation of global preference index, calculation of positive and negative outranking flows, and ranking.

Step1: Determination of deviations based on pairwise comparisons (Equation 12)

$$d_j(a,b) = g_j(a) - g_j(b) \quad (12)$$

Where $d_j(a,b)$ is the difference between the evaluations of *a* and *b* on each criterion

Step 2: Application and selection of the preference function (Equation 13)

$$P_j(a,b) = F_j[d_j(a,b)] \, j = 1,..k. \quad (13)$$

Where $P_j(a,b)$ denotes the preference of alternative *a* with regard to alternative *b* on each criterion as function of $d_j(a,b)$

Step 3: Calculation of an overall or global preference index (Equation 14)

$$\forall \, a,b \in A, \pi(a,b) = \sum_{j=1}^{k} P_j(a,b) \, w_j \quad (14)$$

Where $\pi(a,b)$ of a over b is defined as the weighted sum $P_j(a,b)$ of for each criterion, and $w_j$ is the weight associated with $j^{th}$ criterion

Step 4: Calculation of outranking flow/the PROMOTHEE Partial Ranking Equations (15 & 16)

$$\varphi^+(a) = \sum_{x \in A} \pi(a,b) \quad (15)$$
$$\varphi^-(a) = \sum_{x \in A} \pi(b,a) \quad (16)$$

Where $\varphi^+(a)$ and $\varphi^-(a)$ denote the positive outranking and negative outranking flow for each alternative, respectively

Step 5: Calculation of net outranking flow/ the PROMETHEE complete ranking (Equation 17)

$$\varphi(a) = \varphi^+(a) - \varphi^-(a) \qquad (17)$$

Where $\varphi(a)$ denote the net outranking flow for each alternative.

Step 6: Select the best/suitable alternative having the highest net value

PROMETHEE can be easily implemented due to its computational efficiency in prioritising jobs in a large-scale problem. Practising engineers can implement this approach with ease to make quick decisions. The evaluation matrix lists all the production parameters for developing the priority rule. In almost all the practical scenarios, not all parameters would be of the same trend, i.e. some parameters need to be minimized, whereas the others need to be maximized. A good solution may be a balance of all the parameters under consideration based on the business requirements. Since all data are in various trends, it is always essential to normalize the available data to obtain valuable and meaningful interpretations.

Table 5: Decision matrix for MK1 instances

| Criterion weights | 0.10 | 0.26 | 0.15 | 0.04 | 0.45 |
|---|---|---|---|---|---|
| Jobs | $C_r1$ | $C_r2$ | $C_r3$ | $C_r4$ | $C_r5$ |
| J 1 | 27.0 | 49.5 | 6 | 1.3 | 3.8 |
| J 2 | 20.0 | 37.5 | 5 | 1.2 | 3.5 |
| J 3 | 27.0 | 48.0 | 5 | 1.6 | 4.2 |
| J 4 | 22.0 | 40.5 | 5 | 1.3 | 3.7 |
| J 5 | 34.0 | 60.0 | 6 | 1.7 | 4.3 |
| J 6 | 26.0 | 45.0 | 6 | 1.3 | 3.2 |
| J 7 | 17.0 | 33.0 | 5 | 1.2 | 3.2 |
| J 8 | 33.0 | 49.5 | 5 | 1.6 | 3.3 |
| J 9 | 24.0 | 45.0 | 6 | 1.2 | 3.5 |
| J 10 | 25.0 | 46.5 | 6 | 1.2 | 3.6 |

Note: **$C_r1$**- Processing Time, **$C_r2$**- Due Date, **$C_r3$**- Number of operations, **$C_r4$**-Setup time, **$C_r5$**- Slack Time per Remaining Per Operation;

The weights assigned for the criteria based on their relative importance in the ranking have been presented. Compute each job's preference index and net φ value for every job. The outranked indices of jobs among the alternatives and their corresponding ranks are shown for the MK1 problem instance, for a numerical illustration of PROMETHEE is given in Appendix 2.

## 4.2 Decision Matrix and Job Ranking of all MCDM Methods

An example of the four MCDM methods employed for ranking a benchmark problem, MK1, instances (10 jobs and 6 machines) is illustrated in this section. The evaluation matrix is shown in Table 5. The weights of the Criteria determine the magnitude of their effect. The evaluation matrix provides the inputs for all the MCDMs. Initial ranks for a benchmark instance MK1 using TOPSIS, EDAS, CP and PROMETHEE are presented in Table 6. The ranking of jobs is dynamically changing throughout the simulation. An external spreadsheet is used for computing the job priorities and interfaced with the DES model. The input for the MCDM method includes criteria selected for ranking and their weights calculated using FAHP. The ranking of jobs is performed based on the logic explained in the pseudocode.

Table 6: Initial Job Ranking for TOPSIS, EDAS, CP and PROMETHEE

| Jobs/ MCDM | TOPSIS | EDAS | CP | PROMETHEE |
|---|---|---|---|---|
| J 1 | 9 | 9 | 8 | 10 |
| J 2 | 2 | 7 | 4 | 4 |
| J 3 | 8 | 10 | 10 | 3 |
| J 4 | 3 | 4 | 7 | 2 |
| J 5 | 10 | 3 | 9 | 9 |
| J 6 | 4 | 5 | 2 | 6 |
| J 7 | 1 | 1 | 1 | 1 |
| J 8 | 7 | 2 | 3 | 8 |
| J 9 | 5 | 8 | 5 | 6 |
| J 10 | 6 | 6 | 6 | 5 |

## *4.3 Integration of MCDMs with DES*

Integration of DES with MCDM is required to evaluate the performance of PFJSS instances. Initially, weights were assigned to the criteria using FAHP based on the chosen criteria. This acts as the input for ranking jobs using MCDM, which serves as the input for DES. The schema of integration of the simulation model with MCDM is shown in Figure 4.

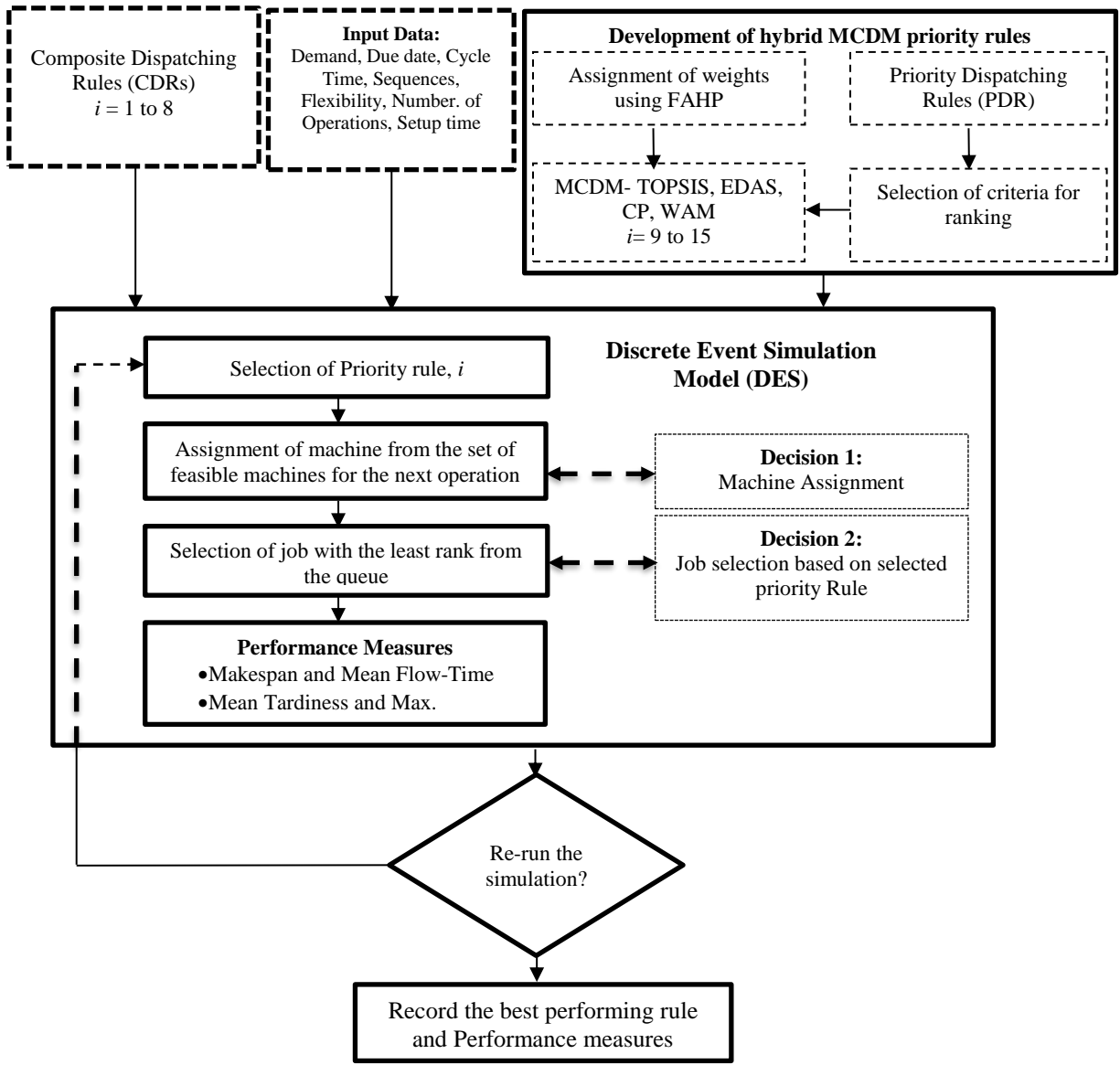

Figure 4 Caption: Integration of MCDM and Simulation Model

Figure 4 Alt Text: Figure depicting the communication happening between the Simulation model and the MCDM, which is the characteristic feature of the study

Table 7: Comparative assessment of measures for various benchmark instances with different arrival patterns

| Arrival Pattern | | Instances/Measures | MK1 | MK2 | MK3 | MK4 | MK5 | MK6 | MK7 | MK8 | MK9 | MK10 |
|---|---|---|---|---|---|---|---|---|---|---|---|---|
| Static Arrival Pattern | | Makespan | 52.35$^{C9}$ | 52.13$^{C12}$ | 256.46$^{C5}$ | 85.61$^{C1}$ | 199.13$^{C5}$ | 112.04$^{C12}$ | 214.23$^{C7}$ | 577.43$^{C3}$ | 410.33$^{C7}$ | 324.15$^{C10}$ |
| | | Mean Flow time | 33.46$^{C7}$ | 35.45$^{C2}$ | 165.45$^{C8}$ | 41.06$^{C12}$ | 109.45$^{C9}$ | 95.45$^{C1}$ | 135.45$^{C7}$ | 334.03$^{C5}$ | 262.14$^{C11}$ | 230.60$^{C10}$ |
| | | Mean Tardiness | 8.46$^{C5}$ | 14.45$^{C7}$ | 45.54$^{C12}$ | 18.64$^{C5}$ | 41.70$^{C9}$ | 23.45$^{C5}$ | 56.45$^{C2}$ | 65.06$^{C11}$ | 76.54$^{C10}$ | 65.16$^{C12}$ |
| | | Max Tardiness | 20.58$^{C12}$ | 28.51$^{C5}$ | 120.5$^{C10}$ | 26.13$^{C9}$ | 76.09$^{C11}$ | 25.62$^{C7}$ | 136.92$^{C12}$ | 125.64$^{C6}$ | 195.07$^{C9}$ | 160.32$^{C2}$ |
| Dynamic Arrival pattern | Uniform | Makespan | 47.45$^{C5}$ | 48.56$^{C7}$ | 240.14$^{C1}$ | 60.57$^{C2}$ | 180.67$^{C7}$ | 105.45$^{C9}$ | 182.15$^{C10}$ | 574.05$^{C12}$ | 357.90$^{C15}$ | 300.49$^{C11}$ |
| | | Mean Flow time | 27.41$^{C9}$ | 35.64$^{C1}$ | 150.45$^{C5}$ | 34.65$^{C7}$ | 95.83$^{C12}$ | 87.01$^{C9}$ | 123.45$^{C5}$ | 316.55$^{C12}$ | 243.50$^{C8}$ | 214.73$^{C11}$ |
| | | Mean Tardiness | 5.53$^{C5}$ | 6.26$^{C2}$ | 17.63$^{C12}$ | 5.14$^{C10}$ | 10.37$^{C9}$ | 8.26$^{C12}$ | 15.95$^{C10}$ | 22.20$^{C11}$ | 32.03$^{C6}$ | 26.95$^{C8}$ |
| | | Max Tardiness | 8.42$^{C3}$ | 13.04$^{C11}$ | 71.56$^{C10}$ | 13.56$^{C5}$ | 28.51$^{C7}$ | 25.65$^{C12}$ | 32.36$^{C6}$ | 43.56$^{C8}$ | 113.49$^{C11}$ | 60.35$^{C5}$ |
| | Increasing | Makespan | 51.45$^{C7}$ | 48.91$^{C1}$ | 232.64$^{C5}$ | 75.40$^{C11}$ | 178.35$^{C7}$ | 110.94$^{C9}$ | 180.05$^{C12}$ | 580.30$^{C4}$ | 384.08$^{C9}$ | 307.36$^{C12}$ |
| | | Mean Flow time | 30.77$^{C11}$ | 35.18$^{C8}$ | 154.65$^{C2}$ | 35.47$^{C7}$ | 96.47$^{C10}$ | 101.61$^{C9}$ | 116.54$^{C5}$ | 312.10$^{C10}$ | 257.65$^{C7}$ | 210.14$^{C12}$ |
| | | Mean Tardiness | 5.27$^{C3}$ | 6.22$^{C6}$ | 38.82$^{C11}$ | 6.53$^{C1}$ | 10.36$^{C9}$ | 12.45$^{C12}$ | 16.45$^{C10}$ | 27.45$^{C5}$ | 35.64$^{C8}$ | 33.25$^{C7}$ |
| | | Max Tardiness | 10.24$^{C10}$ | 15.45$^{C1}$ | 100.36$^{C8}$ | 13.39$^{C5}$ | 25.64$^{C7}$ | 14.36$^{C6}$ | 39.30$^{C9}$ | 45.65$^{C5}$ | 124.54$^{C11}$ | 90.54$^{C8}$ |
| | Decreasing | Makespan | 45.72$^{C5}$ | 47.29$^{C7}$ | 242.87$^{C8}$ | 63.70$^{C9}$ | 186.50$^{C12}$ | 105.63$^{C10}$ | 182.87$^{C5}$ | 530.54$^{C7}$ | 362.71$^{C12}$ | 281.51$^{C9}$ |
| | | Mean Flow time | 34.83$^{C4}$ | 37.34$^{C8}$ | 162.32$^{C2}$ | 35.30$^{C10}$ | 93.89$^{C3}$ | 100.19$^{C6}$ | 116.52$^{C9}$ | 312.08$^{C7}$ | 250.39$^{C12}$ | 215.86$^{C9}$ |
| | | Mean Tardiness | 5.05$^{C11}$ | 6.13$^{C4}$ | 30.57$^{C8}$ | 5.01$^{C4}$ | 10.02$^{C9}$ | 10.21$^{C7}$ | 21.26$^{C5}$ | 30.18$^{C1}$ | 40.43$^{C12}$ | 29.74$^{C10}$ |
| | | Max Tardiness | 9.51$^{C6}$ | 14.66$^{C7}$ | 86.11$^{C10}$ | 11.65$^{C1}$ | 25.17$^{C12}$ | 36.44$^{C9}$ | 36.01$^{C11}$ | 45.21$^{C8}$ | 121.98$^{C12}$ | 95.02$^{C9}$ |
| | Random | Makespan | 62.70$^{C2}$ | 66.46$^{C5}$ | 290.93$^{C1}$ | 94.67$^{C7}$ | 228.77$^{C9}$ | 150.43$^{C2}$ | 224.72$^{C8}$ | 622.74$^{C10}$ | 430.21$^{C7}$ | 359.26$^{C12}$ |
| | | Mean Flow time | 42.12$^{C2}$ | 44.81$^{C8}$ | 171.90$^{C10}$ | 50.09$^{C9}$ | 119.45$^{C7}$ | 119.65$^{C3}$ | 123.16$^{C12}$ | 350.34$^{C6}$ | 276.55$^{C9}$ | 236.78$^{C5}$ |
| | | Mean Tardiness | 8.45$^{C5}$ | 8.57$^{C12}$ | 46.53$^{C3}$ | 10.96$^{C11}$ | 14.35$^{C8}$ | 20.67$^{C10}$ | 26.49$^{C12}$ | 44.50$^{C5}$ | 51.05$^{C7}$ | 42.84$^{C9}$ |
| | | Max Tardiness | 17.49$^{C10}$ | 25.64$^{C4}$ | 122.65$^{C1}$ | 25.45$^{C7}$ | 50.42$^{C9}$ | 50.88$^{C12}$ | 70.98$^{C11}$ | 80.87$^{C8}$ | 142.98$^{C5}$ | 135.98$^{C8}$ |

Note: **C1**-2PT+LWR+FDD, **C2**-2PT+LWR+Slack, **C3**- SPT+LWR+Slack, **C4**-2PT+LWR+EDD, **C5**-(7*LTWC) +(11*SPT)+12*(LNOP+AT), **C6**- LTWC/(3+LNOP-LRNOP), **C7**- ODD+RT, **C8**- [EDD+[(LRNOP+LTWC)/(LRWC-LTWC)]*LNOP]*LRNOP, **C9**-TOPSIS, **C10**- EDAS, **C11**-CP **C12**-PROMETHEE, **PT**- Processing Time, **LWR**- Least Work Remaining, **FDD** – Flow Due Date, **EDD**- Earliest Due Date, **LTWC**- Least Total work content, **SPT**- Shortest Processing Time, **LNOP**- Least No of Operations, **LRWC**- Least Remaining work content, **EDD**- Earliest Due Date, **LRNOP**- Least Remaining no of operations, **PT** – Process Time, **RT**- Remaining Time, **ODD** – Operational Due Date

## 5. Results and Discussion

The performance of the proposed rules is compared with the best-known CDRs to justify the effectiveness of the proposed approach. With this view, the proposed MCDM-based rules and the chosen CDRs from the literature were tested on ten benchmark instances and real-time problem instances for static and different job arrival patterns. The details of the performance evaluation of static and dynamic job arrival scenarios of PFJS benchmark instances are discussed in sections 5.1 and 5.2. The best-performing rule for every problem instance for each performance criterion is shown in Table 7.

### *5.1 Performance Evaluation of static arrival pattern*

With respect to Makespan, the CDR C5: (7*LTWC) +(11*SPT) +12*(LNOP+AT), CDR C7: ODD + Re, and MCDM C12: PROMETHEE produced the best results for two instances each. The problems MK3 & MK5, MK7 & MK9, and MK2 & MK6 are the respective pairs of instances where CDRs C5, C7, and MCDM C12 outperformed the other rules considered in this study. Considering the Mean Flow time, the CDR C7: ODD + RT outmatched the other rules for MK1 and MK7 whilst, for all other instances, different rules were producing the best result for each instance. The rules that fetched the least Mean Flow time for one instance each include CDRs C1, C2, C5, & C8, and MCDMs C9, C10, C11, & C12. Considering the tardiness related measures, C5: (7*LTWC) +(11*SPT) +12*(LNOP+AT) achieved the best results across three instances MK1, MK4, and MK6 while the rules C2, C7, C9, C10, C11, C12 performed better than other rules for one instance each. C9: TOPSIS and C12: PROMETHEE has produced the best result concerning maximum Tardiness for two problem instances, each MK4 & MK9 and MK1 & MK7, respectively. The rules C2, C5, C6, C7, C10, and C11 gave the best values with respect to maximum Tardiness for one instance each. It is observed that for most of the problem instances considering static arrival of jobs for benchmark problems, CDRs are performing well compared to MCDM-based rules. However, CDRs and MCDM-based rules prove to be toe-to-toe overall.

*5.2 Performance Evaluation of dynamic arrival pattern*

CDR C7: ODD + RT has produced the best results for MK2 & MK5, MK1 & MK5, MK2 & MK8, and MK4 & MK9 problem instances for uniform, increasing arrival rate and decreasing arrival rate of arrival, and random arrival patterns, respectively. CDR C5: (7*LTWC) +(11*SPT) +12*(LNOP+AT) produced the best results for MK1 & MK9 for uniform arrival, and MK1 & MK7 for decreasing rate of arrival. MCDM C9: TOPSIS has produced the best results for *MK*6 &*MK*9 for an increasing rate of arrival and *MK*3 &*MK*10 for decreasing rate of arrival. MCDM C12: PROMETHEE has also outperformed other rules for two instances, MK7 and MK10, for the increasing rate of arrival of jobs. For the Mean Flow time, MCDM C9: TOPSIS has produced the best results for dynamic job arrivals in the greatest number of instances. The TOPSIS MCDM rule has given the best results for MK1 & MK6 for uniform job arrivals, MK3 & MK10 for decreasing the rate of job arrivals, and MK4 & MK9 for random job arrivals. Besides, the rule CDR C5: (7*LTWC) +(11*SPT) +12*(LNOP+AT) has produced the least value of Mean Flow time for two instances, MK3 & MK7 for uniform job arrival and MK7 & MK10 for random arrivals. Apart from these, the rules CDR C7: ODD + RT and C10: EDAS has produced the best values for MK4 & MK9 and MK5 & MK8, respectively, for increasing job arrivals. The rules CDR: C2, C4, C5, C6, & C8, and MCDM: C10 and C12 have given the best results for one instance each for decreasing the rate of job arrivals.

The tardiness-based measures show a much larger variation in the rules that perform best for each instance. MCDM: C12: PROMETHEE, CDR: C7: ODD + RT, MCDM: C10: EDAS, CDR: C4: 2PT+LWR+EDD, and CDR: C5: (7*LTWC) +(11*SPT) +12*(LNOP+AT) have performed the best for two instances each with respect to the mean Tardiness. The MCDM rule C12 has given the best results for four instances overall. The MCDM rule produced the best results for MK3 & MK6 for uniform arrival and MK2 & MK7 for random arrival of jobs. MK4 & MK9, MK5 & MK8, MK2 & MK4, and MK1 & MK8 are the respective pairs of instances for which the rules CDR: C7, MCDM: C10, CDR: C4, and CDR: C5 have shown the best values with respect to mean Tardiness.

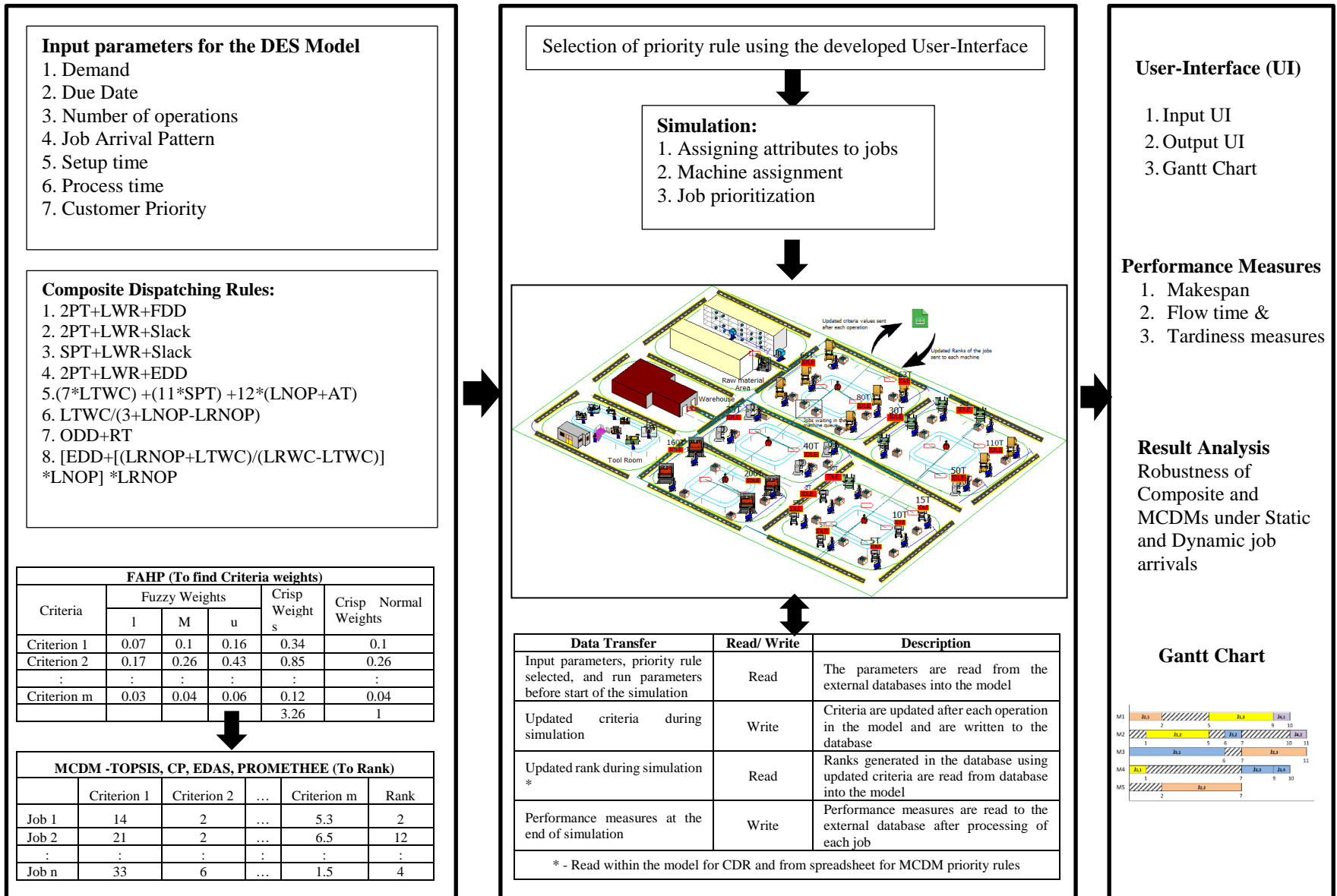

Figure 5 Caption: Work Flow of Real-time Simulation
Figure 5 Alt Text: The real-time simulation framework is presented in three broad steps: data input, arriving job priorities using proposed MCDMs and simulation, and recording the performance measures.

For maximum Tardiness measure, CDR: C2, C5, & C6, and MCDM: C9, C10, & C11 produced best result for uniform job arrival; CDR: C2, C4, C5, C6, & C8, and MCDM: C10 & C12 for decreasing rate of job arrivals; and CDR: C3, C7, & C8, and MCDM: C9, C10, & C11 for random job arrivals. CDR: C5: (7*LTWC) +(11*SPT) +12*(LNOP+AT) and CDR: C8: [EDD+[(LRNOP+LTWC)/(LRWC-LTWC)] *LNOP] *LRNOP showed the best results of maximum Tardiness for four instances each. CDR: C5 showed the best results for MK4 & MK10 and MK4 & MK8 for uniform arrival and increasing arrival rate, respectively. The CDR: C8 produced the best results for the MK3 & MK10 and MK8 and MK10, respectively, for increasing and random arrival of jobs. Apart from these, MCDM: C11: CP, MCDM: C9: TOPSIS, and MCDM: C12: PROMETHEE gave the best values for the instances MK2 & MK9, MK6 & MK10, and MK5 & MK9 respectively. The rules CDR: C1, C4, C5, & C7, and MCDM: C9, C10, C11, and C12 performed the best for one instance each for random arrival scenarios.

The performance of CDRs is better than MCDMs for small-size benchmarking problems. For Makespan, MCDMs perform better than CDRs for 4/10 and 17/40 instances for static and dynamic job arrival patterns, respectively. MCDMs perform better for 4/10 and 18/40 instances for static and dynamic job arrival patterns for Mean Flow Time. Considering mean Tardiness, CDRs and MCDMs produce the best results for five instances each for static arrival patterns. For the dynamic arrivals scenario, it is observed that MCDMs are performing well for 19/40 instances. With respect to the maximum tardiness measure, the proposed MCDMs are achieving the best results for 6/10 and 17/40 instances with static and dynamic arrivals patterns, respectively. Mean and average percentage deviation statistics are computed and observed that MCDMs are consistent in their performance compared to CDRs.

## 6 Case Study

An auto-ancillary press shop has been undertaken to evaluate the effectiveness of the proposed rules. Data relevant to process mapping, machine details, demand patterns, job sequence, job setting time, and job due date are collected from the shop floor. The press shop has 28 workstations (Press) of varying capacities. Each machine's job material, geometry, and tonnage determine the set of jobs that a machine can handle. Demands are considered monthly, and job prioritization is done through past experiences and observations. The high number of machines and jobs increases the complexity and vulnerability of errors in scheduling. A simulation model is developed using DES software to represent the press shop.

Table 8: Comparative assessment of MCDM based PDR and CDR for all performance measures

| Job Arrival Pattern | | Performance Measures/Priority Rules | Composite Dispatching Rules (Time in Hours) | | | | | | | | MCDM Methods (Time in Hours) | | | |
|---|---|---|---|---|---|---|---|---|---|---|---|---|---|---|
| | | | | | | | | | | | TOPSIS | EDAS | CP | PROMETHEE |
| | | | C1 | C2 | C3 | C4 | C5 | C6 | C7 | C8 | C9 | C10 | C11 | C12 |
| Static | | Makespan | 196.35 | 209.66 | 212.38 | 206.05 | 191.50 | 217.13 | 207.54 | 197.82 | **174.65** | 180.54 | 185.34 | 178.68 |
| | | Mean Flow time | 77.38 | 87.16 | 84.98 | 89.09 | 78.03 | 86.06 | 78.33 | 87.22 | **70.14** | 76.57 | 80.32 | 75.64 |
| | | Mean. Tardiness | 26.53 | 21.68 | 21.06 | 23.11 | 22.15 | 23.87 | 26.46 | 26.08 | 20.47 | 18.87 | **13.92** | 16.84 |
| | | Max. Tardiness | 56.24 | 54.69 | 51.63 | 51.79 | 57.76 | 54.46 | 50.42 | 57.53 | 45.34 | 42.62 | **41.34** | 48.47 |
| Dynamic | Uniform | Makespan | 186.68 | 194.8 | 193.15 | 180.12 | 192.47 | 204.47 | 201.38 | 193.64 | 175.45 | 184.62 | 180.14 | **172.64** |
| | | Mean Flow time | 90.64 | 89.21 | 86.47 | 87.19 | 83.36 | 84.21 | 82.08 | 84.20 | **65.47** | 72.89 | 76.17 | 69.72 |
| | | Mean. Tardiness | 28.40 | 32.13 | 26.67 | 24.65 | 25.74 | 26.56 | 24.18 | 21.13 | 19.33 | 17.41 | 14.22 | **12.39** |
| | | Max. Tardiness | 51.09 | 50.50 | 49.39 | 47.39 | 54.27 | 58.43 | 52.06 | 50.19 | 42.82 | 44.17 | **38.47** | 40.61 |
| | Increasing | Makespan | 192.75 | 184.05 | 192.53 | 194.53 | 190.68 | 178.25 | 181.75 | 179.74 | **162.47** | 170.64 | 172.92 | 168.25 |
| | | Mean Flow time | 83.02 | 74.35 | 76.42 | 81.98 | 82.19 | 84.02 | 82.57 | 78.05 | **70.72** | 78.98 | 80.41 | 74.56 |
| | | Mean. Tardiness | 28.34 | 18.75 | 26.53 | 18.11 | 22.38 | 27.73 | 25.73 | 27.71 | 16.28 | **12.47** | 14.45 | 18.63 |
| | | Max. Tardiness | 51.01 | 42.91 | 50.72 | 53.43 | 44.12 | 45.13 | 52.21 | 42.60 | 36.42 | **30.44** | 34.78 | 38.65 |
| | Decreasing | Makespan | 193.83 | 185.20 | 196.23 | 183.79 | 197.29 | 180.85 | 182.82 | 194.56 | **164.52** | 172.62 | 174.22 | 168.47 |
| | | Mean Flow time | 78.58 | 78.40 | 79.00 | 80.28 | 83.67 | 82.14 | 81.41 | 82.43 | 68.12 | 73.45 | 78.62 | **65.32** |
| | | Mean. Tardiness | 30.75 | 28.96 | 25.04 | 29.31 | 23.45 | 26.48 | 23.36 | 29.78 | 20.45 | 17.23 | **14.24** | 18.64 |
| | | Max. Tardiness | 51.38 | 44.55 | 54.83 | 51.41 | 56.17 | 45.03 | 55.72 | 45.94 | 42.15 | 39.64 | **35.41** | 40.27 |
| | Random | Makespan | 211.75 | 222.23 | 230.47 | 219.97 | 214.67 | 232.30 | 224.44 | 226.74 | **182.34** | 197.52 | 192.62 | 188.64 |
| | | Mean Flow time | 95.38 | 92.62 | 90.58 | 99.99 | 97.11 | 98.25 | 99.49 | 94.90 | **76.61** | 81.38 | 85.94 | 79.47 |
| | | Mean. Tardiness | 29.48 | 27.97 | 27.27 | 32.00 | 32.79 | 22.56 | 33.12 | 27.82 | 24.25 | 26.88 | **15.41** | 19.68 |
| | | Max. Tardiness | 66.93 | 56.49 | 65.61 | 68.54 | 62.87 | 55.67 | 56.09 | 65.73 | 50.21 | 54.54 | **42.14** | 48.41 |

Note: **C1**-2PT+LWR+FDD, **C2**-2PT+LWR+Slack, **C3**- SPT+LWR+Slack, **C4**-2PT+LWR+EDD, **C5**-(7*LTWC) +(11*SPT) +12*(LNOP+AT), **C6**-LTWC/(3+LNOP-LRNOP), **C7**- ODD+RT, **C8**- [EDD+[(LRNOP+LTWC)/(LRWC-LTWC)] *LNOP] *LRNOP, **C9**-TOPSIS, **C10**- EDAS, **C11**-CP **C12**-PROMETHEE

The values of processing times are directly provided as input to the simulation model using an external database. Demand and setup times of jobs are generated using the probability distribution function. The scalability of the simulation model in terms of the number of jobs and its parameters makes it applicable to several problems. The real-time simulation framework is schematically represented in Figure 5. It explains three steps involved: data input, arriving job priorities using proposed MCDMs and simulation, and recording the performance measures.

## 6.1 Performance evaluation of static arrival pattern

Makespan recorded the best value of 174.7 hours for TOPSIS, followed by 178.7 hours for PROMETHEE, as shown in Table 8. The CDR: C5: (7*LTWC) +(11*SPT) +12*(LNOP+AT) reported the least Makespan of 191.5 hours. The proposed MCDM rules performed well with respect to Makespan. Mean Flow Time recorded the least value of 70.14 hours for TOPSIS, followed by 75.54 hours for PROMETHEE. Among the CDRs, C1: 2PT+LWR+FDD reported the best value of 77.38 hours for Makespan. Production delays hinder job completion within specified due dates and reduce the efficiency of job shops. Mean Tardiness had the least value of 13.92 hours for MCDM: CP, with 16.84 hours being the second least value for PROMETHEE. The minimum value of Mean Tardiness for CDRs was reported to be 21.06 hours for C3- SPT+LWR+Slack. MCDM: CP was the best-performing hybrid rule with respect to Maximum Tardiness with a value of 41.34 hours followed by MCDM: EDAS with 42.62 hours. Concerning Maximum Tardiness, CDR: C7 ODD+RT produced the best result of 50.42 hours. MCDMs rule outmatched the CDRs across all performance measures for the chosen case study. The proposed hybrid PDRs effectively improve all selected performance measures involving a number of operations, due date, process time, demand rate, setup time, and customer priority.

## 6.2 Performance evaluation of dynamic arrival pattern

For the four different dynamic arrival patterns used, the performance measures were recorded to evaluate the effectiveness and robustness of the patterns. TOPSIS produced the least Makespan values of 162.47 hrs for the increasing rate of arrivals followed by 164.52 hrs for decreasing rate of arrival pattern. Mean Flow Time was the minimum for decreasing arrival pattern with a value of 65.32 hrs for the PROMETHEE rule. The second-best value of 65.47 hrs was achieved for TOPSIS for a uniform arrival pattern. Mean Tardiness produced the least value of 12.39 hrs for PROMETHEE when a uniform arrival pattern was followed. Increasing arrival pattern gave a closer value of 12.47 hrs for the EDAS priority rule. The fourth measure, maximum Tardiness, showed the best value of 30.44 hrs and 35.44 hrs for increasing and decreasing arrival patterns, respectively, for the rules EDAS and CP.

It is observed that the MCDM-based rules have outperformed the Composite rules for all the dynamic arrival patterns with respect to all performance measures. This is because the MCDM rules have been developed considering multiple parameters that affect scheduling. This also offers the flexibility to alter the parameters and their weights. Also, this simulation-based approach that accounts for varying job arrivals is more helpful for the scheduling managers to incorporate any changes in the arrivals and make quick decisions based on the simulation results. The determination of the robustness of the composite rules and MCDMs for the static and dynamic arrivals followed this.

The robustness of each method is summarized with respect to the various performance measure. Since TOPSIS involves the normalization of criteria to bring them to comparable scales and it performs well for Makespan and mean Flow Time measures, the role of AHP is significant in assigning proper weights to these criteria. In the present study, customer priority has been given the highest weightage, followed by the number of operations and processing time. This indicates the strong dependence of Makespan and Mean Flow Time on these three criteria. Since customer priority has the utmost importance, followed by a number of operations and due date, and EDAS produces the best results for Tardiness-based measures, it could be attributed that these criteria are vital in determining Tardiness-based measures. Since the ranking of jobs is being done without normalizing the criteria in CP, the criterion such as demand and due date that are relatively higher in magnitude than the other criteria are being given more weightage for the ranking of jobs. Given that CP performs better for Tardiness-based measures, it could be inferred that demand and due date play a vital role in these measures. PROMETHEE considers the minimum and maximum values of each criterion which are used to create the preference function by comparing them with the criterion values of each instance. This rule has produced the best results for one instance each in Makespan and mean Flow Time and also mean Tardiness.

*6.3 Effects of criteria weights on performance measures*

In the present study, the weights for the criteria have been assigned using Fuzzy Analytic Hierarchy Process (FAHP). In the Flexible Job Shop, various factors that influence the performance of the system are identified based on surveys and discussions with industry experts. A questionnaire was prepared and distributed to the shop floor engineers and managers for feedback on factors influencing the shop's performance. Delphi technique (Dalkey and Helmer 1963) has been used to conduct the survey and identify factors affecting shop performance. The criteria chosen in the present study are customer priority, demand, due date, setup time, number of operations, processing time, and slack time remaining per operation. From the scores obtained for each of the parameters using Delphi, ranks were assigned in descending order. The weights / relative importance of criteria is obtained using the FAHP process.

The performance of the shop floor was evaluated by assigning different weights for criteria. Weights are chosen in such a way that their sum is equal to 1. Since many experiments are required to assess the shop floor's performance, experiments were carried out by varying the weight of one factor at a time, starting from 0.1 to 0.9 in steps of 0.1. The design table for conducting the experiments is given in Table 1A (Appendix 3). This procedure was carried out considering the MCDM-TOPSIS method for the case problem considered in this study. The results are presented in Figure 6.

When the weights of process time increased, Makespan and mean Flow Time were reduced, whereas the Tardiness-based measures were increased. Contrarily, when the weights of due date were increased, Makespan and Mean Flow Time were increased while the Tardiness-based measures decreased. It is noticed that increasing the weights for Setup time, number of operations and customer priority does not influence the output performance measures. However, there were slight improvements in maximum Tardiness with increased weights of a number of operations and setup time. An increase in weights of demand improved the performance with respect to Makespan, Mean Flow time and Maximum Tardiness. There was a slight increase in mean Tardiness with the increase in the weight of demand. The present study's results match the trends available in the literature regarding the influence of weights over performance measures. (Holthaus and Rajendran,1997).

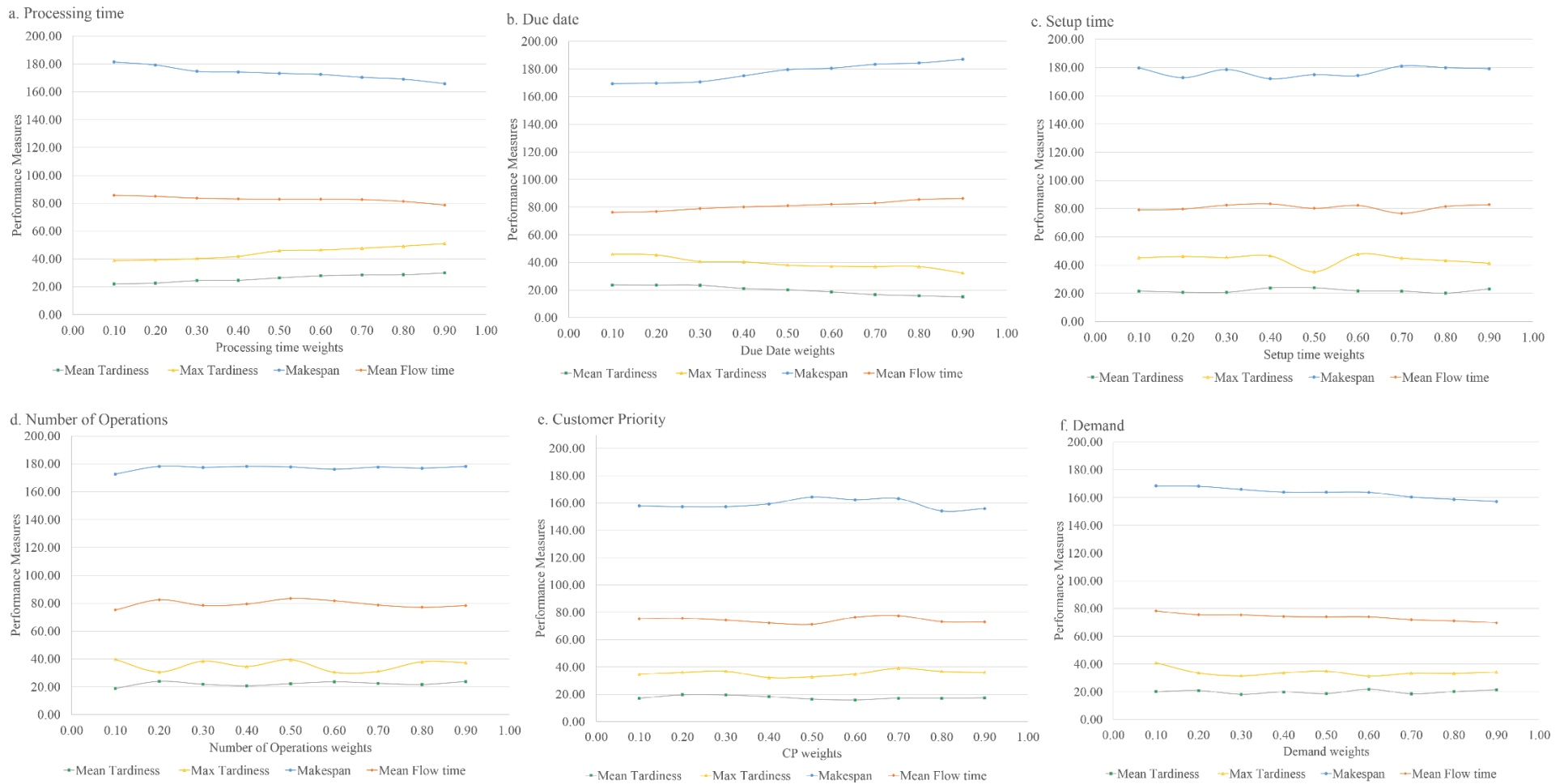

Figure 6 Caption: Effects of criterion weights on performance measures a) process time, b) Due date, c) Setup time, d) Number of Operations, e) Customer Priority, f) Demand

Figure 6 Alt Text: Graphs representing the variations in the measured performance indicators with respect to changes in the criterion weights

*6.4 Statistical Validation*

To compare the differences in performance of the CDRs and MCDM rules, two statistical tests, 'ANOVA' and 't-test', were performed. To carry out these tests, the CDRs were considered as a group, whereas each MCDM was considered as a separate entity. The ANOVA test (Table 9) reveals that CDRs and MCDMs are significantly different with respect to performance measures considered in this study. It is observed that for all performance measures, the p-value value computed is less than 0.05.

Table 9: Result of one-way ANOVA Test

|  |  | Sum of Squares | df | Mean Square | F | Sig. |
|---|---|---|---|---|---|---|
| Makespan | Between Groups | 6793.986 | 4 | 1698.497 | 9.132 | 9.90E-06 |
|  | Within Groups | 10229.827 | 55 | 185.997 |  |  |
|  | Total | 17023.813 | 59 |  |  |  |
| Mean_Flow_Time | Between Groups | 1700.185 | 4 | 425.046 | 11.550 | 7.08E-07 |
|  | Within Groups | 2024.088 | 55 | 36.802 |  |  |
|  | Total | 3724.272 | 59 |  |  |  |
| Mean_Tardiness | Between Groups | 1060.718 | 4 | 265.180 | 20.930 | 1.57E-10 |
|  | Within Groups | 696.832 | 55 | 12.670 |  |  |
|  | Total | 1757.550 | 59 |  |  |  |
| Max_Tardiness | Between Groups | 1838.233 | 4 | 459.558 | 11.753 | 5.73E-07 |
|  | Within Groups | 2150.544 | 55 | 39.101 |  |  |
|  | Total | 3988.777 | 59 |  |  |  |

Note: **df** – Degrees of freedom, **F**- Test Statistic, **Sig.**-P value

The Mean values obtained for each performance measure through independent sample t-tests (Table 10) show that TOPSIS had the least mean score for Makespan and Mean Flow time while CP had the least mean scores for the tardiness-related measures. It is evident that the MCDMs have outperformed the CDRs for all the performance measures. The p-values of TOPSIS, EDAS and PROMETHEE are less than 0.05, which shows a significant difference between the use of CDRs and other MCDMs with respect to each outcome measure. Whereas the p-value for Mean Flow Time corresponding to CP is 0.108, which is greater than 0.05, indicating that there is no significant difference between the performances of CDR and CP with respect to Mean Flow time (Table 11)

Table 10: Mean Performance measures

| Individual_MCDM | **CDRs** | **TOPSIS** | **EDAS** | **CP** | **PROMETHEE** |
|---|---|---|---|---|---|
| Makespan | 199.4140 | 171.8860 | 181.1880 | 181.0480 | 175.3360 |
| Mean_Flow_Time | 85.3110 | 70.2120 | 76.6540 | 80.2920 | 72.9420 |
| Mean_Tardiness | 26.1455 | 20.1560 | 18.5720 | 14.4480 | 17.2360 |
| Max_Tardiness | 53.3233 | 43.3880 | 42.2820 | 38.4280 | 43.2820 |

Table11: Independent sample t-test

| t-test | | TOPSIS | | | EDAS | | | CP | | | PROMETHEE | | |
|---|---|---|---|---|---|---|---|---|---|---|---|---|---|
| | | Levene's Test for Equality of Variances | | t-test for Equality of Means | Levene's Test for Equality of Variances | | t-test for Equality of Means | Levene's Test for Equality of Variances | | t-test for Equality of Means | Levene's Test for Equality of Variances | | t-test for Equality of Means |
| | | F | Sig. | Sig. (2-tailed) | F | Sig. | Sig. (2-tailed) | F | Sig. | Sig. (2-tailed) | F | Sig. | Sig. (2-tailed) |
| Makespan | Equal variances assumed | 2.340 | 0.133 | 0.000 | 1.421 | 0.240 | 0.013 | 2.630 | 0.112 | 0.011 | 2.363 | 0.132 | 0.001 |
| | Equal variances not assumed | | | 0.000 | | | 0.014 | | | 0.003 | | | 0.001 |
| Mean_Flow_Time | Equal variances assumed | 1.985 | 0.166 | 0.000 | 1.950 | 0.170 | 0.007 | 2.758 | 0.104 | 0.108 | 0.294 | 0.590 | 0.000 |
| | Equal variances not assumed | | | 0.000 | | | 0.002 | | | 0.031 | | | 0.004 |
| Mean_Tardiness | Equal variances assumed | 0.938 | 0.338 | 0.001 | 0.254 | 0.617 | 0.000 | 6.297 | 0.016 | 0.000 | 0.586 | 0.448 | 0.000 |
| | Equal variances not assumed | | | 0.006 | | | 0.030 | | | 0.000 | | | 0.001 |
| Max_Tardiness | Equal variances assumed | 0.551 | 0.462 | 0.002 | 0.189 | 0.666 | 0.001 | 1.475 | 0.231 | 0.000 | 0.183 | 0.671 | 0.002 |
| | Equal variances not assumed | | | 0.007 | | | 0.044 | | | 0.000 | | | 0.005 |

## 7 Conclusions

DES and MCDM approaches have been integrated into this study to determine job priorities in large-scale PFJS problems involving multiple criteria. Benchmark instances of various sizes and a real-world PFJS with 114 job varieties and 28 machines were used to examine the performance of the proposed approach. This work aims to rank the jobs by considering real-world criteria such as slack time, processing time, number of operations, demand, customer priority, setup time, and due date. Criteria are identified using the Delphi method, and weights for the criteria are obtained using FAHP. MCDMs are used to rank the jobs whilst the DES model is used to evaluate the performance measures considering static and dynamic job arrival patterns. Makespan, Flow Time, Mean Tardiness, and Maximum Tardiness measures are used to assess the shop floor performance. Proposed methods are compared with the best-performing CDRs available in the literature for Makespan, Flow time, and Tardiness measures. For benchmark problem instances, CDRs outperformed the MCDMs for smaller-size problems. For Makespan, MCDMs performed better than CDRs for 40% and 42.5% of problem instances for static and dynamic job arrival patterns, respectively. Similar results were observed for mean Flow Time, where MCDMs performed better for 40% and 45% of problem instances for static and dynamic job arrival patterns, respectively. Considering Mean Tardiness, CDRs and MCDMs produce the best results for half of the problem instances each for static arrival patterns. For the dynamic arrivals scenario, it is observed that MCDMs are performing well for 47.5% of problem instances. With respect to Maximum Tardiness, the MCDMs outmatched CDRs for 60% and 42.5% of problem instances with static and dynamic arrivals patterns, respectively. For real-world problem instances, MCDMs performed well compared to CDRs for static and dynamic arrival conditions. For static arrival, C9: TOPSIS and C11: Compromise Programming produced the best results for Flow Time and Tardiness measures. In the case of dynamic job arrival patterns, C9: TOPSIS and C12: PROMETHEE achieved the best results for Makespan and Flow time measures. C12: PROMETHEE produced better results for mean Tardiness, whereas C10: EDAS and C11: Compromise programming produced better results for both measures of Tardiness.

Overall it is observed that, for a large-scale instance, the MCDM approach outperforms CDRs for the performance measures considered in this study. The better performance of MCDM methods can be attributed to the ranking of jobs based on demand rate, due date, customer priority, setup time, process time, and the number of operations that influence the scheduling in the real world. The DES model developed in this study dynamically updates the parameters and job ranking as the simulation progresses. The methodology proposed in this study is modular and simple for implementation in an industrial scenario. The proposed model is adaptive to incorporate additional parameters such as the number of resources, job variants, static and dynamic job arrival patterns, and flexibility. The DES model provides information about the job priorities, performance measures, and utilization of resources based on the proposed MCDMs. This would help the industries and policymakers to identify a pattern in terms of the parameters that influence the scheduling effectiveness with respect to the problem instances studied, thereby helping them make decisions concerning the selection of suitable MCDM. The model proposed in this study is generic and can be used to analyze scenarios involving probabilistic activity times and demand patterns. The tailor-made MCDM rules that function effectively based on the competencies of the PDRs using which they have been built would be of more interest to the researchers and policymakers as they encompass the inclusion of all the criteria required to be considered for scheduling. Further, developing a decision support system might have more practical significance for production managers to schedule jobs effectively in a more complex environment.

**Data availability statement**

The data that support the findings of this study are available from the corresponding author upon reasonable request.

# Appendix 1

## A) Pseudocode to generate job ranking using TOPSIS method

The pseudocode of the TOPSIS method used for prioritizing jobs is given in this section.

Input:

$x_{ij}$ = value of criterion '$j$' of job instance '$i$'// Job instances

$w_j$ = weight of criterion '$j$'//Criterion's weight computed using FAHP

Output:

$R_i$ // Ranks obtained for chosen job instances

Makespan, Mean Flow Time, Mean Tardiness, and Max Tardiness // Performance Measures

**start**
**while** ($x_{ij} \ne$ NULL) **do**
**for** i= 1 to *n* **do** //Create a Decision Matrix *D*
    **for** $x_{ij}$ in *D* **do** $r_{ij} = \dfrac{x_{ij}}{\sqrt{\sum_{i=1}^{m} X_{ij}^2}}$ // Compute normalized decision matrix
    **end for**
    $V_{ij} = w_j \cdot r_{ij}$ // Weighted Normalized decision matrix, $w_j$ is weights of each criterion calculated by FAHP
    **if** $j \epsilon J$ **then**
        $v_i^+ = \{max(v_{ij}) \; if \; j\epsilon J; min(v_{ij})$
        else if $j\epsilon J$ then
        $v_i^- = \{min(v_{ij}) \; if \; j\epsilon J; max(v_{ij})$
    **end if**
    **for** $J \epsilon \; j \; \epsilon \; C_j$ **do**
        $S_i^+ = \sqrt{[\sum(v_i^* - v_{ij})^2]}$ // Calculate the separation measure of PIS ($S_i^+$)
        $S_i^- = \sqrt{[\sum(v_i^- - v_{ij})^2]}$// Calculate the separation measure of NIS ($S_i^-$)
    **end for**
    **for** each $x_{ij}$ in *D* do
        $cc_i = \dfrac{S_i^+}{(S_i^+ + S_i^-)}$ //Calculate the relative closeness coefficient
    **end for**
    Rank(*n*) based on *cci* // Rank the jobs based on closeness coefficient
**end for**
    **do**
    Simulation
    Assigning ranks to jobs in the DES model
    Updating ranks of jobs
**end while**
Evaluation of the performance measures
Generation of Gantt chart
**Stop**

### B) Pseudocode to generate job ranking using EDAS method

The pseudocode of the EDAS method used for prioritizing jobs is given in this section.

Input:
 $x_{ij}$ = value of criterion '$j$' of job instance '$i$' // Job instances
 $w_j$ = weight of criterion '$j$' //Criterion's weight computed using FAHP

Output:
 $R_i$ // Ranks obtained for chosen job instances
 Makespan, Mean Flow Time, Mean Tardiness, and Max Tardiness // Performance Measures

**start**
**while** ($x_{ij} \neq$ NULL) **do**
**for** i= 1 to *n* **do** //Create a Decision Matrix *D*
  **for** $x_{ij}$ in *AV* **do** $AVJ = \frac{\sum_{i=1}^{n} X_{ij}}{n}$ // Computing the Average Solution Matrix
  **end for**
  **for** $x_{ij}$ in *AV* **do**
   $PDA_{ij} = \frac{\max(0,(Xij-AVJ))}{AVJ}$, $NDA_{ij} = \frac{\max(0,(AVJ-Xij))}{AVJ}$ // Computing PDA and NDA
  **end for**
  **for** $x_{ij}$ in *PDA* **do**
   $SP_i = \sum_{j=1}^{m} wjPDAij$ // Computation of the weighted sum of positive ($SP_i$)
  **end for**
  **for** $x_{ij}$ in *NDA* **do**
   $SN_i = \sum_{j=1}^{m} wjNDAij$ // Computation of the weighted sum of negative ($SN_i$)
  **end for**
  **for** $x_{ij}$ in *NSP* **do**
   $NSP_i = \frac{SPi}{max(SPi)}$ // Normalizing the weighted sum of *PDA*
  **end for**
  **for** $x_{ij}$ in *NSN* **do**
   $NSN_i = 1 - \frac{SPi}{max(SPi)}$ // Normalizing the weighted sum of *NDA*
  **end for**
  **for** $x_{ij}$ in $AS_i$ **do**
   $AS_i = (NSP_i - NSN_i) \times 0.5$ // Computing the Appraisal Score (AS)
  Rank(*n*) based on $AS_i$ // Rank the jobs based on AS
  **end for**
**end for**
 **do**
 Simulation// DES model
 Assigning ranks to jobs in the DES model
 Updating ranks of jobs
**end while**
Evaluation of the performance measures
Generation of Gantt chart
**stop**

## C) Pseudocode to generate job ranking using CP method

The methodology used for developing job ranking using CP is given in this section. Job ranking is computed dynamically at the end of completion of every operation. The pseudocode for arriving the job ranking is given in the form of pseudocode in this section.

Pseudocode:

Input:
    $x_{ij}$ = value of criterion '$j$' of job instance '$i$' // Job instances
    $w_j$ = weight of criterion '$j$' //Criterion's weight computed using FAHP

Output:
    $R_i$ // Ranks obtained for chosen job instances
    Makespan, Mean Flow Time, Mean Tardiness, and Max Tardiness // Performance Measures

**start**
**while** ($x_{ij} \neq$ NULL) **do**
**for** $i = 1$ to $n$ **do** //Create a Decision Matrix $D$
    **for** all $C_i$ ($j = 1$ to n) // Calculate the normalized weights

$$L_P(a) = \left[ \sum_{j=1}^{J} w_j \left| \frac{f_j^* - f_j(a)}{f_j^* - f_j^{**}} \right|^P \right]^{\frac{1}{P}}$$ 
// Compute the '$Lp$' Matrix

    Rank($n$) based on $L_p$ // Rank the jobs based on the least metric value
    **end for**
**end for**
    **do**
    Simulation// DES model
    Assigning ranks to jobs in the DES model
    Updating ranks of jobs
**end while**
Evaluation of the performance measures
Generation of Gantt chart
**Stop**

*D) Pseudocode to generate job ranking using PROMETHEE method*

The Pseudocode of only the PROMETHEE method used for prioritizing jobs is given in this section. Other Pseudocode of TOPSI, EDAS and CP is given in the Appendix 1.

**Input:**
$x_{ij}$ = value of criterion '$j$' of job instance '$i$'// Job instances
$w_j$ = weight of criterion '$j$'//Criterion's weight computed using FAHP

**Output:**
$R_i$ // Ranks obtained for chosen job instances
Makespan, Mean Flow Time, Mean Tardiness, and Max Tardiness// Performance Measures

**start**
**while** ($X_{ij} \neq$ NULL) **do**
**for** i= 1 to *n* **do** //Create a Decision Matrix *D*
    $d_j(a,b) = g_j(a) - g_j(b)$ // Determination of deviations based on pairwise comparisons
    **for** all *j* = 1 to *n do*
        $P_j(a,b) = F_j[d_j(a,b)]$ $(j = 1,..k)$ // Computing of the preference function
        $\forall\ a, b \in A,\ \pi(a,b) = \sum_{j=1}^{k} P_j(a,b)\ w_j$ // Calculation of global preference index
    **end for**
    **for** all $x \in A$ **do**
    $\varphi^+(a) = \sum_{x \in A} \pi(a,b)$ *and* $\varphi^-(a) = \sum_{x \in A} \pi(b,a)$ // Calculation of outranking flow/the PROMETHEE II partial ranking
        $\varphi(a) = \varphi^+(a) - \varphi^-(a)$// Calculation of net outranking flow/ the PROMETHEE II complete ranking
    **end for**
    Rank(*n*) based on $\varphi(a)$ // Rank the jobs based on the highest net value
**end for**
    **do**
    Simulation// DES model
    Assigning ranks to jobs in the DES model
    Updating ranks of jobs
**end while**
Evaluation of the performance measures

**stop**

# Appendix 2

*Numerical illustration of PROMETHEE method*

The approach to calculating the net φ and ranking the alternatives is presented. Table 2A depicts the decision matrix for calculating the rank of each job using PROMETHEE

Table 2A: Decision Matrix of PROMETHEE for MK1 instance

| Job/Criteria | Process time | Due Date | Number of operations | Setup time | STROP |
|---|---|---|---|---|---|
| 1 | 27 | 50 | 6 | 1.4 | 3.8 |
| 2 | 20 | 38 | 5 | 1.2 | 3.5 |
| 3 | 27 | 48 | 5 | 1.6 | 4.2 |
| 4 | 22 | 41 | 5 | 1.3 | 3.7 |
| 5 | 34 | 60 | 6 | 1.7 | 4.3 |
| 6 | 26 | 45 | 6 | 1.4 | 3.2 |
| 7 | 17 | 33 | 5 | 1.0 | 3.2 |
| 8 | 33 | 50 | 5 | 1.7 | 3.3 |
| 9 | 24 | 45 | 6 | 1.2 | 3.5 |
| 10 | 25 | 47 | 6 | 1.3 | 3.6 |

Note: J1-J10- Job set of MK1 instances (Brandimarte, 1993), STROP- Slack Time Remaining Per Operation

Table 2B: Transformed payoff matrix of PROMETHEE

| Job/Criteria | Process time | Due Date | Number of operations | Setup time | STROP |
|---|---|---|---|---|---|
| J1 | -27 | -50 | -6 | -1.4 | 3.8 |
| J2 | -20 | -38 | -5 | -1.2 | 3.5 |
| J3 | -27 | -48 | -5 | -1.6 | 4.2 |
| J4 | -22 | -41 | -5 | -1.3 | 3.7 |
| J5 | -34 | -60 | -6 | -1.7 | 4.3 |
| J6 | -26 | -45 | -6 | -1.4 | 3.2 |
| J7 | -17 | -33 | -5 | -1.0 | 3.2 |
| J8 | -33 | -50 | -5 | -1.7 | 3.3 |
| J9 | -24 | -45 | -6 | -1.2 | 3.5 |
| J10 | -25 | -47 | -6 | -1.3 | 3.6 |

Note: J1-J10- Job set of MK1 instances (Brandimarte, 1993), STROP- Slack Time Remaining Per Operation

To allow the analysis of the problem from a maximization viewpoint, a negative sign is assigned to the minimization criteria as (-min) = max. The payoff matrix post-transformation of criteria to maximization type is presented in Table 2B.

**Step 1:** Pairwise difference between values of alternatives for criteria

The foremost step in the computation of rank is a pairwise comparison of criteria. For instance, for criterion $C_1$, the pairwise difference between alternative $J_1$ and $J_2$ are -27-(-20) = -7. Likewise, the pairwise difference between alternatives $J_2$ and $J_1$ for $C_1$ is -20-(-27) =7. Also, the difference between $J_1$ and $J_1$ would be 0 as the comparison is made between the same criteria. Pairwise differences between the alternatives for criteria $C_1$ to $C_5$ are calculated.

**Step 2:** The values of the preference function for $C_1$ to $C_5$ were based on the quasi-criterion function. An illustration based on $J_1$ and $J_2$ is discussed. Pairwise differences between alternatives $J_1$ and $J_2$ for criterion $C_1$ are – 27 – (-20) = -7. For quasi criterion function with a threshold value of 0, preference function value $PF_1(J_1, J_2)$ =0 (as -7 <20). Similarly difference for alternatives $J_2$ and $J_1$ are -20- (-27) = 7 and corresponding preference function value $PF_1(J_2, J_1)$ =1 (as 7 <0).

**Step 3**: Computation of multi-criterion preference index

Multi-criterion preference index (Table 2C), $\pi(J_1, J_2)$ for pairwise alternative $(J_1, J_2)$ is calculated as follows.

Preference function values for $J_1$ and $J_2$ for criteria $C_1$ to $C_5$ are 0,0,1.0,0.5 and 0, respectively. The corresponding weights of the criteria are 0.1, 0.26, 0.15, 0.04 and 0.45

$$\pi(J_1 J_2) \frac{(7X0X0.1) + (5X0.5X0.26) + (0X0X0.15) + (0.35X0.5X0.04) + (0.5X0.5X0.45)}{0.1 + 0.26 + 0.15 + 0.04 + 0.45} = 1.83$$

Similarly,

$$\pi(J_2 J_1) \frac{(-7X0.5X0.1) + (2X0.5X0.26) + (1X1X0.15) + (0.15X0.5X0.04) + (0.3X0.5X0.45)}{0.1 + 0.26 + 0.15 + 0.04 + 0.45} = 2.13$$

Table 2C: Multi-criterion Preference Index values

|  | J1 | J2 | J3 | J4 | J5 | J6 | J7 | J8 | J9 | J10 |
|---|---|---|---|---|---|---|---|---|---|---|
| J1 | 0.00 | 1.83 | 0.10 | 0.00 | 1.24 | 0.00 | 0.00 | 0.31 | 0.00 | 0.00 |
| J2 | 2.13 | 0.00 | 1.88 | 0.54 | 0.00 | 1.43 | 1.22 | 2.22 | 1.33 | 1.59 |
| J3 | 0.35 | 2.08 | 0.00 | 0.00 | 0.78 | 0.15 | 0.00 | 0.50 | 0.15 | 0.15 |
| J4 | 1.59 | 3.43 | 1.34 | 0.00 | 0.00 | 0.94 | 1.43 | 1.73 | 0.84 | 1.08 |
| J5 | 0.00 | 3.87 | 0.45 | 0.00 | 1.24 | 0.00 | 0.00 | 0.00 | 0.00 | 0.00 |
| J6 | 0.77 | 2.60 | 0.67 | 0.11 | 0.07 | 0.00 | 0.54 | 0.96 | 0.07 | 0.29 |
| J7 | 2.94 | 4.77 | 2.69 | 1.34 | 0.81 | 2.17 | 0.00 | 2.98 | 2.13 | 2.40 |
| J8 | 0.26 | 1.79 | 0.20 | 0.09 | 0.05 | 0.15 | 0.62 | 0.00 | 0.20 | 0.22 |
| J9 | 0.81 | 2.64 | 0.71 | 0.05 | 0.00 | 0.10 | 0.00 | 1.04 | 0.00 | 0.27 |
| J10 | 0.54 | 2.37 | 0.44 | 0.02 | 0.00 | 0.05 | 0.00 | 0.80 | 0.00 | 0.00 |

Note: J1-J10- Job set of MK1 instances (Brandimarte, 1993)

Table 2D: Job Priority Ranking of the PROMETHEE method

| S.No. | Φ+ (a,b) | Φ- (a,b) | Net(Φ) | Rank |
|---|---|---|---|---|
| 1 | 2.24 | 9.38 | -7.14 | 8 |
| 2 | 15.08 | 0.92 | 14.16 | 2 |
| 3 | 3.38 | 8.02 | -4.65 | 7 |
| 4 | 10.94 | 2.15 | 8.79 | 3 |
| 5 | 0.00 | 25.49 | -25.49 | 10 |
| 6 | 5.54 | 4.99 | 0.56 | 5 |
| 7 | 22.22 | 0.00 | 22.22 | 1 |
| 8 | 3.02 | 10.54 | -7.52 | 9 |
| 9 | 6.42 | 4.70 | 1.72 | 4 |
| 10 | 4.27 | 5.99 | -1.73 | 6 |

Note: Φ+ (a,b) – Positive outranking flow, Φ- (a,b)- Negative outranking flow, Net(Φ)-Net outranking

**Step 4** Computation of φ+ (as per equation)

$\varphi+(J_1) = 0 + 1.83 + 0.1 + 0 + 0 + 0.31 + 0 + 0 = 2.24$

**Step 5** Computation of φ+ (as per equation)

$\varphi-(J_1) = 0 + 2.13 + 0.35 + 1.59 + 0.77 + 0 + 2.94 + 0.26 + 0.81 + 0.54 = 9.38$

**Step 6**: Computation of net φ (as per equation)

Netφ $(J1) = \varphi+(J_1) - \varphi-(J_1) = 2.24 - 9.38 = -7.14$

Similarly, other net φ values were calculated and ranked the jobs. From Table 2D, it is evident that $J_7$ is considered as rank 1 as with the highest net value of 22.22 is considered as the best. The procedure illustrated in the numerical example has been implemented to determine job priorities of all benchmark problem instances considering static and dynamic job arrival patterns. DES models are used to determine the performance measures considered in this study to evaluate the proposed priority rule.

*Table1A: Design of Experiments*

| Scenarios | **PT** | DD | ST | NOP | CP | D | Scenarios | PT | **DD** | ST | NOP | CP | D |
|---|---|---|---|---|---|---|---|---|---|---|---|---|---|
| 1 | **0.10** | 0.18 | 0.18 | 0.18 | 0.18 | 0.18 | 1 | 0.18 | **0.10** | 0.18 | 0.18 | 0.18 | 0.18 |
| 2 | **0.20** | 0.16 | 0.16 | 0.16 | 0.16 | 0.16 | 2 | 0.16 | **0.20** | 0.16 | 0.16 | 0.16 | 0.16 |
| 3 | **0.30** | 0.14 | 0.14 | 0.14 | 0.14 | 0.14 | 3 | 0.14 | **0.30** | 0.14 | 0.14 | 0.14 | 0.14 |
| 4 | **0.40** | 0.12 | 0.12 | 0.12 | 0.12 | 0.12 | 4 | 0.12 | **0.40** | 0.12 | 0.12 | 0.12 | 0.12 |
| 5 | **0.50** | 0.10 | 0.10 | 0.10 | 0.10 | 0.10 | 5 | 0.10 | **0.50** | 0.10 | 0.10 | 0.10 | 0.10 |
| 6 | **0.60** | 0.08 | 0.08 | 0.08 | 0.08 | 0.08 | 6 | 0.08 | **0.60** | 0.08 | 0.08 | 0.08 | 0.08 |
| 7 | **0.70** | 0.06 | 0.06 | 0.06 | 0.06 | 0.06 | 7 | 0.06 | **0.70** | 0.06 | 0.06 | 0.06 | 0.06 |
| 8 | **0.80** | 0.04 | 0.04 | 0.04 | 0.04 | 0.04 | 8 | 0.04 | **0.80** | 0.04 | 0.04 | 0.04 | 0.04 |
| 9 | **0.90** | 0.02 | 0.02 | 0.02 | 0.02 | 0.02 | 9 | 0.02 | **0.90** | 0.02 | 0.02 | 0.02 | 0.02 |
| Scenarios | PT | DD | **S.T** | NOP | CP | D | Scenarios | P.T | DD | ST | **NOP** | CP | D |
| 1 | 0.18 | 0.18 | **0.10** | 0.18 | 0.18 | 0.18 | 1 | 0.18 | 0.18 | 0.18 | **0.10** | 0.18 | 0.18 |
| 2 | 0.16 | 0.16 | **0.20** | 0.16 | 0.16 | 0.16 | 2 | 0.16 | 0.16 | 0.16 | **0.20** | 0.16 | 0.16 |
| 3 | 0.14 | 0.14 | **0.30** | 0.14 | 0.14 | 0.14 | 3 | 0.14 | 0.14 | 0.14 | **0.30** | 0.14 | 0.14 |
| 4 | 0.12 | 0.12 | **0.40** | 0.12 | 0.12 | 0.12 | 4 | 0.12 | 0.12 | 0.12 | **0.40** | 0.12 | 0.12 |
| 5 | 0.10 | 0.10 | **0.50** | 0.10 | 0.10 | 0.10 | 5 | 0.10 | 0.10 | 0.10 | **0.50** | 0.10 | 0.10 |
| 6 | 0.08 | 0.08 | **0.60** | 0.08 | 0.08 | 0.08 | 6 | 0.08 | 0.08 | 0.08 | **0.60** | 0.08 | 0.08 |
| 7 | 0.06 | 0.06 | **0.70** | 0.06 | 0.06 | 0.06 | 7 | 0.06 | 0.06 | 0.06 | **0.70** | 0.06 | 0.06 |
| 8 | 0.04 | 0.04 | **0.80** | 0.04 | 0.04 | 0.04 | 8 | 0.04 | 0.04 | 0.04 | **0.80** | 0.04 | 0.04 |
| 9 | 0.02 | 0.02 | **0.90** | 0.02 | 0.02 | 0.02 | 9 | 0.02 | 0.02 | 0.02 | **0.90** | 0.02 | 0.02 |
| Scenarios | PT | DD | ST | NOP | **CP** | D | Scenarios | PT | DD | ST | NOP | CP | **D** |
| 1 | 0.18 | 0.18 | 0.18 | 0.18 | **0.10** | 0.18 | 1 | 0.18 | 0.18 | 0.18 | 0.18 | 0.18 | **0.10** |
| 2 | 0.16 | 0.16 | 0.16 | 0.16 | **0.20** | 0.16 | 2 | 0.16 | 0.16 | 0.16 | 0.16 | 0.16 | **0.20** |
| 3 | 0.14 | 0.14 | 0.14 | 0.14 | **0.30** | 0.14 | 3 | 0.14 | 0.14 | 0.14 | 0.14 | 0.14 | **0.30** |
| 4 | 0.12 | 0.12 | 0.12 | 0.12 | **0.40** | 0.12 | 4 | 0.12 | 0.12 | 0.12 | 0.12 | 0.12 | **0.40** |
| 5 | 0.10 | 0.10 | 0.10 | 0.10 | **0.50** | 0.10 | 5 | 0.10 | 0.10 | 0.10 | 0.10 | 0.10 | **0.50** |
| 6 | 0.08 | 0.08 | 0.08 | 0.08 | **0.60** | 0.08 | 6 | 0.08 | 0.08 | 0.08 | 0.08 | 0.08 | **0.60** |
| 7 | 0.06 | 0.06 | 0.06 | 0.06 | **0.70** | 0.06 | 7 | 0.06 | 0.06 | 0.06 | 0.06 | 0.06 | **0.70** |
| 8 | 0.04 | 0.04 | 0.04 | 0.04 | **0.80** | 0.04 | 8 | 0.04 | 0.04 | 0.04 | 0.04 | 0.04 | **0.80** |
| 9 | 0.02 | 0.02 | 0.02 | 0.02 | **0.90** | 0.02 | 9 | 0.02 | 0.02 | 0.02 | 0.02 | 0.02 | **0.90** |

*Note: **PT**-Processing Time, **DD**- Due Date; **ST**-Setup Time; **NOP**-NO of Operations; **CP**-Customer priority; **D**-Demand*